\documentclass[a4paper,leqno]{article}

\usepackage{epsf}

\usepackage[dvips]{graphicx}
\usepackage{amsmath}
\usepackage{amssymb}
\usepackage{amsfonts}
\usepackage{amsthm}
\usepackage{enumerate}

\theoremstyle{plain}
\newtheorem{teo}{Theorem}[section]
\newtheorem{lem}[teo]{Lemma}

\newtheorem{prop}[teo]{Proposition}

\theoremstyle{definition}
\newtheorem{defin}[teo]{Definition}
\newtheorem{oss}[teo]{Remark}
\newtheorem{exam}[teo]{Example}
\newtheorem*{assA}{Assumption A}

\renewcommand{\eqref}[1]{\textnormal{(\ref{#1})}}

\numberwithin{equation}{section}

\newcommand{\cvd}{\hfill$\square$}
\renewcommand{\proof}[1]{\noindent\textsc{Proof#1}}

\newcommand{\rmd}{\mathrm{d}}

\title{Reconstruction of material losses by perimeter penalization and phase-field methods}

\author{Luca \textsc{Rondi}\\
\normalsize{Universit\`a degli Studi di Trieste}\\
\normalsize{Dipartimento di Matematica e Informatica}\\
\normalsize{via Valerio, 12/1}\\
\normalsize{34127 Trieste, Italy}\\
\normalsize{\texttt{rondi@units.it}}}

\date{}

\begin{document}

\maketitle

\setcounter{section}{0}
\setcounter{secnumdepth}{1}

\begin{abstract}
We treat the inverse problem of determining material losses, such as cavities, in a conducting body, by performing electrostatic measurements at the boundary. We develop a numerical approach, based on variational methods, to reconstruct the unknown material loss by a single boundary measurement of current and voltage type.

The method is based on the use of phase-field functions to model the material losses and on a perimeter-like penalization to regularize the otherwise ill-posed problem.
We justify the proposed approach by a convergence result, as the error on the measurement goes to zero.

\medskip


\medskip

\end{abstract}

\section{Introduction}

In many inverse or optimal shape problems arising in the applications, the aim is to reconstruct the shape of an object, usually represented by an unknown open set, satisfying certain requirements. If we restrict ourselves to a variational formulation, for the sake of simplicity, we look for the shape minimizing a given functional $F$ among all the admissible shapes. The shape is often modeled as a binary function, that is the open set is described through its characteristic function.

Two of the main issues for a satisfactory numerical resolution of this kind of problems are the following. First of all, and especially for inverse problems, the problem may be ill-posed, that is stability is missing or, in other words, $F$ is not continuous.
Second, numerically handling shapes or sets is not an easy task from the implementation point of view.
The first issue is usually tackled by a regularization method, namely by adding to the functional a term penalizing the binary function with respect to some $BV$-related norm. For most applications, this should be enough for ensuring a regularization without being restricting on the class of admissible unknowns. Often the $BV$-related norm is simply a perimeter-like penalization.
About the handling of shapes or sets in computations, in many cases this is performed by associating to the open set a smooth function describing it. For example, one way of doing it is to replace the characteristic function of an open set $D$ with a smooth function, referred to as a \emph{phase-field function},
which is close to $0$ outside $D$, close to $1$ inside $D$, and has a quick transition from $0$ to $1$ across the boundary of $D$. Another way is the so-called level-set method, where $D$ is identified with the sublevel set $\{\psi<0\}$ of a smooth function
$\psi$.

We are interested in using perimeter-like regularizations and phase-field functions
for solving inverse or optimal shape problems, in particular those that are not well-posed. We aim to prove in a rigorous way that this kind of approach provides a good approximation of the original problem, allowing us at the same time to tame the ill-posedness and to have a formulation amenable to be easily implemented.
A cornerstone of this method is the approximation, in the sense of $\Gamma$-convergence, of the perimeter functional by functionals defined on phase-field functions, due to Modica and Mortola, \cite{Mod e Mor}. Since \cite{Mod} such a result has found innumerable applications. In fact, whenever the functional $F$ is continuous in a suitable way, the invariance of $\Gamma$-convergence by continuous perturbations permits to obtain an analogous $\Gamma$-convergence result if we add to $F$ the perimeter penalization. Whenever the problem is ill-posed, that is $F$ is not continuous,
a corresponding convergence result is not straightforward any more. Since we believe
that the method is valuable also in the ill-posed case, it would be important to justify it in a rigorous way, in general through a convergence result inspired by $\Gamma$-convergence techniques, for various interesting applications.

In this paper we perform such an analysis for the following inverse problem, arising from non-destructive evaluation. We aim to
determine perfectly insulating defects in a homogeneous and isotropic conducting body by performing electrostatic measurements of voltage and current type at the boundary.
The conducting body is contained in 
$\Omega$,
a bounded
domain of $\mathbb{R}^N$, $N\geq 2$. The defects may have different geometrical properties,
for instance we may have at the same time
\emph{cracks} (either interior or surface breaking), or \emph{material losses} (either interior, that is cavities, or at the boundary).
We denote with $K$ the union of the boundaries of
these defects, whereas $\tilde{\gamma}$ is  
a part of the boundary of $\Omega$ which is accessible,
known and disjoint from $K$. If a current density
$f\in L^2(\tilde{\gamma})$, with zero mean,
is applied on $\tilde{\gamma}$, then
the electrostatic potential $u=u(f,K)$
is the solution to the following Neumann boundary value problem
\begin{equation}\label{dirpbm0}
\left\{\begin{array}{ll}
\Delta u=0&\text{in }\Omega\backslash K,\\
\nabla u\cdot \nu=f&\text{on }\tilde{\gamma},\\
\nabla u\cdot \nu=0&
\text{on }\partial(\Omega\backslash K)\backslash\tilde{\gamma}.
\end{array}
\right.
\end{equation}
We call $G_K$ the connected component of $\Omega\backslash K$ 
which is reachable from $\tilde{\gamma}$ and
we say that a defect is a \emph{material loss} if $G_K$ is equal
to the interior of its closure, that is if no crack-type defect is present.

The value of $u$, that is the voltage, may be measured on another 
part of the boundary of $\Omega$, say $\gamma$, which we assume to be accessible,
known, disjoint from $K$ and belonging to $\partial G_K$.
We call $g$ such a measurement, that is $g=u|_{\gamma}$. For simplicity, we may also assume that $\gamma$ coincides with $\tilde{\gamma}$.
If the defect $K$ is unknown, we aim to recover its shape and location,
that is $G_K$, by
prescribing one or more current densities $f$ and measuring the corresponding voltage
on $\gamma$, $g=u|_{\gamma}$, where $u$ solves \eqref{dirpbm0}. In mathematical words, we are given one or more pairs of Cauchy data $(g,f)$ on a known part of the boundary and we aim to reconstruct the domain of validity of the elliptic equation.

Here we are interested in the reconstruction only of material losses, that is cavities or material losses at the boundary, and for simplicity we refer to it as the
\emph{inverse cavity problem}. It is well-known that, in every dimension, a single boundary measurement is enough to reconstruct a material loss, thus providing uniqueness for the inverse problem, see for instance \cite{Ron06} for a proof with minimal regularity assumptions on the unknown material loss.
Stability results have been proved in \cite{Ale e Ber e Ros e Ves} for the three dimensional case and in \cite{Ale e Ron} for the planar case, where also the instability character of the problem has been explicitly shown.

We notice that $u$, the electrostatic potential solution to \eqref{dirpbm0}, is constant on any connected component of $\Omega\backslash K$ different from $G_K$. The key observation is that its jump set in $\Omega$ is essentially contained in $K$. The uniqueness result recalled before actually allows us to say more, in fact the jump set of $u$ uniquely determines
$G_K$, that is the unknown material loss. Therefore we are interested in the reconstruction of the electrostatic potential $u$ and especially of its (unknown) jump set.
This suggests the possibility to set up a reconstruction procedure by solving a free-discontinuity problem related to the function $u$.

The main difficulty for the reconstruction is due to the ill-posedness of the problem. 
In fact, since they are measured, the Cauchy data that are available are not exact. Since the problem is severely ill-posed, such an error on the measurements may lead to a much greater error on the reconstructed defect. Furthermore, the inverse problem is nonlinear. In fact, even if the direct problem \eqref{dirpbm0} is linear, the dependence of 
the electrostatic potential $u$, and of its values on $\gamma$, from the defect $K$ is nonlinear. Finally, from a numerical point of view, the fact that the unknown is a set, namely $G_K$, introduces an additional complication for the implementation.

We propose a variational method to tackle at the same time these difficulties. The idea is to use a perimeter-like penalization to regularize the problem and to replace
the unknown set $G_K$ with its characteristic function and, in turn, with a phase-field function, to obtain a formulation that may be implemented numerically. Namely, the regularization we propose is related to the so-called Modica-Mortola functional, an approximation of the perimeter when phase-field functions are used.
We might construct a family of functionals, depending on the noise level on the measurements $\varepsilon$, to be minimized with respect to the variable $u$ (the reconstructed potential) and the phase-field variable $v$. However, to simplify the implementation we would rather have a functional depending on the phase-field variable $v$ only. Thus, we take $u$ depending on $v$, $u=u(v)$, as a solution to an almost degenerate elliptic problem whose coefficient is given by a slight modification of $v$, depending on $\varepsilon$. In other words, we replace the direct cavity problem with an elliptic problem in $\Omega$ where the coefficient of the equation is close to $1$ in $G_K$, close to $0$ outside $G_K$, with a quick transition across the boundary $G_K$. The method consists then of minimizing the so-obtained functionals, depending on $\varepsilon$, with respect to the phase-field variable $v$ only. We remark that the reconstructed material loss may be simply computed by a suitable thresholding of the minimizing phase-field and that an approximation of the looked-for electrostatic potential is given by $u=u(v)$ where $v$ is the minimizing phase-field.

The main result of the paper, Theorem~\ref{convteo2}, is that
the corresponding minimizers $v_{\varepsilon}$
converge, as $\varepsilon\to 0^+$, to the characteristic function of $G_K$, thus identifying the looked-for material loss, and that $u_{\varepsilon}=u(v_{\varepsilon})$ 
converge to the looked-for potential $u$. Such a convergence result, whose proofs is obtained by techniques borrowed by $\Gamma$-convergence, provides a
rigorous justification of the method. About the material loss to be reconstructed, this is assumed to satisfy a Lipschitz type regularity. We finally remark that the method makes use of a single measurement and that is enough 
to reconstruct the whole unknown material loss $K$.

If we instead allow the unknown defect not to be a material loss, that is it may include crack-type defects, for simplicity we refer to this problem as the \emph{inverse crack problem}. About uniqueness, stability and reconstruction results on the inverse crack problem,
we refer to \cite{Bry e Vog04} and the references therein. The main difference between the two cases is that for the determination of cracks one measurement is not enough, however, at least in the planar case, two suitably chosen measurements are sufficient.
In \cite{Ron07,Ron08} a corresponding variational approach for the inverse crack problem has been developed. Again such an approach makes use of a penalization on the $(N-1)$-dimensional measure of the defects and of phase-field functions. Namely, it was constructed
a family of functionals, again depending on the noise level on the measurements $\varepsilon$, 
to be minimized with respect to the variable $u$ (the reconstructed potential) and the phase-field variable $v$. Instead of the perimeter functional and the Modica-Mortola functional, the Mumford-Shah functional \cite{Mum e Sha89} and its approximation, in the sense of $\Gamma$ convergence, due to Ambrosio and Tortorelli \cite{Amb e Tor1,Amb e Tor2} were used, respectively. Also in this case a convergence result guaranteed a justification of the method. In Section~\ref{compsec} we recall the results
obtained in \cite{Ron08} for the inverse crack problem and we compare with
those obtained here for the inverse cavity problem. 
The main difficulty in the implementation of the method of \cite{Ron08} is that
the functional to be minimized depends on two variables, the variable $u$, which should approximate the electrostatic potential, and the variable $v$, which is the phase-field variable that should approximate the jump set of the potential and hence the defect.
It would be desirable to formulate the problem depending on one variable only, for instance only on the phase-field variable. Unfortunately such a formulation, which is proved here for the material loss case, may not be feasible.
In fact, Section~\ref{compsec} is devoted to show that the result in \cite{Ron08} is essentially optimal, through several counterexamples.
Moreover, more regularity is needed for the unknown defects of crack-type, namely a regularity assumptions of $C^1$ type, instead of Lipschitz, have to be imposed.
Thus, we show that restricting ourselves to the reconstruction of material losses allows us to gain the following advantages. First, we may lower the a priori assumptions on the unknown defect to ones which are more suited for applications. More importantly,  we obtain and justify a formulation which looks more natural and quite simpler to be implemented.

We finally wish to mention that a numerical implementation, based on the results of this paper and on those of \cite{Ron08}, may be found in \cite{Rin e Ron}. The corresponding
numerical experiments show the validity of these methods also from a practical point of view.

The plan of the paper is the following. After a preliminaries section, Section~\ref{prelsec},
we describe the setting of the direct and inverse problem in
Section~\ref{dirsec}. We treat the material loss case in 
in Section~\ref{cavssec}, where there is the main result of the
paper, Theorem~\ref{convteo2}. In Section~\ref{compsec}, we recall
the results for the inverse crack problem proved in \cite{Ron08} and 
we compare the crack and material loss cases and 
discuss their differences. In particular we show the optimality of the result of \cite{Ron08}.
Finally, in Section~\ref{diffsec} we deal with the differentiability of the functionals involved. Such differentiability is crucial for developing the algorithm used in \cite{Rin e Ron}.

\subsection{Acknowledgments}
This work is partially supported by GNAMPA under 2008 and 2009 projects and by the Italian Ministry of University and Research under PRIN 2008 project. The author wishes also to thank Giuseppe Di Fazio for pointing him out the class of
strong $A_{\infty}$-weights.

\section{Preliminaries}\label{prelsec}

Throughout the paper the integer $N\geq 2$ will denote the space dimension.
 We remark that we shall sometimes
drop the dependence of
any constant upon $N$, the space dimension.
For every $x\in\mathbb{R}^N$, we shall set $x=(x',x_N)$, where $x'\in\mathbb{R}^{N-1}$
and $x_N\in\mathbb{R}$,
and, for any $r>0$, we shall denote by $B_r(x)$ the open ball
in $\mathbb{R}^N$ centred at $x$ of radius $r$.
Usually we shall write $B_r$ instead of $B_r(0)$.
For any subset $E\subset\mathbb{R}^N$ and any $r>0$, we denote $B_r(E)=\bigcup_{x\in E}B_r(x)$.

For any non-negative integer
$k$ we denote by $\mathcal{H}^k$ the $k$-dimensional Hausdorff measure.
For Borel subsets of $\mathbb{R}^N$ the $N$-dimensional
Hausdorff measure coincides with $\mathcal{L}^N$, the $N$-dimensional Lebesgue measure.
Furthermore, if $\gamma\subset\mathbb{R}^N$ is a smooth manifold of dimension $k$,
then $\mathcal{H}^k$ restricted to
$\gamma$ coincides with its $k$-dimensional surface measure.
For any Borel $E\subset\mathbb{R}^N$ we let $|E|=\mathcal{L}^N(E)$.

We recall that a bounded open set $\Omega\subset\mathbb{R}^N$ is said to have a
\emph{Lipschitz boundary}
if for every $x\in\partial\Omega$ there exist a Lipschitz
function $\varphi:\mathbb{R}^{N-1}\to\mathbb{R}$ and a positive constant $r$
such that for any $y\in B_r(x)$ we have, up to a rigid transformation,
$$y\in\Omega\quad \text{if and only if}\quad  y_N<\varphi(y').$$
We observe that in this case the boundary of $\Omega$ has finite $(N-1)$-dimensional Hausdorff measure,
that is $\mathcal{H}^{N-1}(\partial\Omega)<+\infty$.

We say that a function $\varphi:A\to B$, $A$ and $B$ being metric spaces, is \emph{bi-Lipschitz} if it is injective and
$\varphi$ and $\varphi^{-1}:\varphi(A)\to A$ are both Lipschitz functions.
If both the Lipschitz constants of $\varphi$ and $\varphi^{-1}$ are 
bounded by $L\geq 1$, then we say that $\varphi$ is \emph{bi-Lipschitz} with constant $L$.

We recall some basic notation and properties
of functions of bounded variation and sets of finite perimeter. For a more comprehensive treatment of
these subjects see, for instance, \cite{Amb e Fus e Pal, Eva e Gar, Giu}.

Given a bounded open set $\Omega\subset \mathbb{R}^N$,
we denote by $BV(\Omega)$ the Banach space of \emph{functions of bounded
variation}. We recall that $u\in BV(\Omega)$ if and only if  
$u\in L^1(\Omega)$ and its distributional derivative $Du$ is a bounded
vector
measure. We endow $BV(\Omega)$ with the standard norm as follows. Given
$u\in BV(\Omega)$, we denote by $|Du|$ the total variation of its
distributional derivative and
we set $\|u\|_{BV(\Omega)}=\|u\|_{L^1(\Omega)}+|Du|(\Omega)$.
We recall that whenever $u\in W^{1,1}(\Omega)$,
then $u\in BV(\Omega)$ and $|Du|(\Omega)=\int_{\Omega}|\nabla u|$, therefore
$\|u\|_{BV(\Omega)}=\|u\|_{L^1(\Omega)}+\|\nabla u\|_{L^1(\Omega)}=\|u\|_{W^{1,1}(\Omega)}$.

We say that a sequence of $BV(\Omega)$ functions $\{u_h\}_{h=1}^{\infty}$
\emph{weakly}$^*$ \emph{converges} in $BV(\Omega)$ to $u\in BV(\Omega)$ if and only if
$u_h$ converges to $u$ in $L^1(\Omega)$ and $Du_h$
weakly$^*$ converges to $Du$ in $\Omega$, that is
\begin{equation}\label{weakstarconv}
\lim_{h}\int_{\Omega}v \rmd Du_h=
\int_{\Omega}v \rmd Du\quad\text{for any }v\in C_0(\Omega).
\end{equation}
By Proposition~3.13 in \cite{Amb e Fus e Pal}, we have that if
a sequence of $BV(\Omega)$ functions $\{u_h\}_{h=1}^{\infty}$ is bounded in $BV(\Omega)$ and converges to $u$ in $L^1(\Omega)$, then
$u\in BV(\Omega)$ and $u_h$ converges to $u$
weakly$^*$ in $BV(\Omega)$.

Let $\Omega$ be a bounded open set with Lipschitz boundary.
A sequence of $BV(\Omega)$ functions $\{u_h\}_{h=1}^{\infty}$
such that $\sup_h\|u_h\|_{BV(\Omega)}<+\infty$ admits a subsequence
converging weakly$^*$ in $BV(\Omega)$ to a function $u\in BV(\Omega)$, see for instance Theorem~3.23 in \cite{Amb e Fus e Pal}.
As a corollary, we infer that for any $C>0$ the set
$\{u\in BV(\Omega): \|u\|_{BV(\Omega)}\leq C\}$
is a compact subset of $L^1(\Omega)$.

Let $E$ be a bounded 
Borel set contained in $B_R\subset \mathbb{R}^N$. We shall denote by $\chi_E$ its characteristic function. We notice that 
$E$ is compactly
contained in $B_{R+1}$, which we shall denote by $E\Subset B_{R+1}$.
We say that $E$ is a \emph{set of finite perimeter} if
$\chi_E$ belongs to $BV(B_{R+1})$ and we call the number
$P(E)=|D\chi_E|(B_{R+1})$ its \emph{perimeter}.

Let us further remark that the intersection of two sets of finite perimeter is
still a set of finite perimeter. Moreover,
whenever $E$ is open and $\mathcal{H}^{N-1}(\partial E)$ is finite, then $E$ is a set of finite
perimeter, see for instance \cite[Section~5.11, Theorem~1]{Eva e Gar}.
Therefore a bounded open set
$\Omega$ with Lipschitz boundary
is a set of finite perimeter and its perimeter $P(\Omega)$ coincides with
$\mathcal{H}^{N-1}(\partial\Omega)$.

For any bounded open set $\Omega$, we define the following \emph{perimeter functional}
$P:L^1(\Omega)\to [0,+\infty]$ such that
\begin{equation}\label{Pdef}
P(u)=
\left\{
\begin{array}{ll}
\vphantom{\displaystyle{\int}}c|Du|(\Omega)
&\text{if }u\in BV(\Omega)\ \text{and}\ u\in\{0,1\}\text{ a.e.},\\
\vphantom{\displaystyle{\int}}+\infty&\text{otherwise},
\end{array}
\right.
\end{equation}
where $c$ is a positive constant to be chosen later.
We observe that $P(u)=c P(E)$ if $u=\chi_E$ and $E$ is a set of finite perimeter compactly contained in $\Omega$.

We denote by $SBV(\Omega)$ the space of \emph{special functions of bounded
variation}. For any $u\in SBV(\Omega)$, the
density of the absolutely continuous part of $Du$ with respect to $\mathcal{L}^N$ will be denoted
by $\nabla u$, the \emph{approximate gradient} of $u$.
The singular part, with respect to $\mathcal{L}^N$, of $Du$ is
concentrated on $J(u)$, $J(u)$ being the \emph{approximate discontinuity set} (or \emph{jump set}) of $u$ in $\Omega$.
We further say that a function $u\in GSBV(\Omega)$, the space of \emph{generalized
functions of bounded variation}, if $u\in L^1(\Omega)$ and for any $T>0$ its \emph{truncation}
$u_T=(-T)\vee (T\wedge u)\in SBV(\Omega)$. Let us recall that the approximate gradient $\nabla u$ of $u\in GSBV(\Omega)$
is defined almost everywhere and coincides with $\nabla u_T$ almost everywhere on $\{u=u_T\}$, and that
$J(u)=\bigcup_{T>0}J(u_T)$.

The special functions of bounded variation have important compactness
and semicontinuity properties,
see for instance \cite[Theorem~4.7 and Theorem~4.8]{Amb e Fus e Pal}.

We remark that if $u\in BV(\Omega)$ and $u\in\{0,1\}$ almost everywhere in $\Omega$,
then $u\in SBV(\Omega)$ and $P(u)=c |Du|(\Omega)=c \mathcal{H}^{N-1}(J(u))$.

Let us define the so-called \emph{Mumford-Shah functional}, introduced in \cite{Mum e Sha89}
in the context of image segmentation. Let us fix positive constants $b$ and $c$.
Let $\mathcal{MS}:L^1(\Omega)\to [0,+\infty]$ be given by
\begin{equation}\label{MS0}
\mathcal{MS}(u)=\displaystyle{b\int_{\Omega}|\nabla u|^2+c\mathcal{H}^{N-1}(J(u))}\quad
\text{if }u\in GSBV(\Omega),
\end{equation}
whereas $\mathcal{MS}(u)=+\infty$ otherwise.

Let us introduce the following $\Gamma$-convergence results
concerning the approximation of the perimeter functional and the Mumford-Shah functional by phase-field functionals. For the definition and properties of $\Gamma$-convergence we refer to \cite{DaM}. 
The perimeter approximation is due to
Modica and Mortola, \cite{Mod e Mor}, whereas the Mumford-Shah functional approximation is due to Ambrosio and Tortorelli, \cite{Amb e Tor1, Amb e Tor2}.
We shall follow the notation and proofs contained in
\cite{Bra}.

Throughout the paper, for any $p$, $1\leq p\leq +\infty$, we shall denote its conjugate exponent by $p'$, that is $p^{-1}+(p')^{-1}=1$.
Let $W:\mathbb{R}\to[0,+\infty)$ be a continuous function such that
$W(t)=0$ if and only if $t\in\{0,1\}$. Let $c_W=\int_0^1\sqrt{W(t)}\mathrm{d}t$.
In the definition of the perimeter functional we pick $c=2c_W$. For instance, we may choose $W(t)=9t^2(t-1)^2$ for any $t\in \mathbb{R}$, whence $c=2c_W=1$.

The following approximation result is due to Modica and Mortola, \cite{Mod e Mor}, see also \cite[Theorem~4.13]{Bra}.

\begin{teo}\label{Mod-Morteo}
Let $\Omega\subset\mathbb{R}^N$ be a bounded open set with Lipschitz boundary.

For any $\eta>0$ we define the functional
$P_{\eta}:L^1(\Omega)\to [0,+\infty]$ as follows
\begin{equation}\label{modmordef}
P_{\eta}(v)=\left\{
\begin{array}{ll}
\displaystyle{\frac{1}{\eta}\int_{\Omega}W(v)+\eta\int_{\Omega}|\nabla v|^2}&\text{if }v\in W^{1,2}(\Omega,[0,1]),\\
\vphantom{\displaystyle{\int}}+\infty&\text{otherwise}.
\end{array}
\right.
\end{equation}

Then we have that, with respect to the metric of $L^1(\Omega)$,
$P_{\eta}$ $\Gamma$-converges to $P$
as $\eta\to 0^+$.
\end{teo}

Here $W^{1,2}(\Omega,[0,1])=\{v\in W^{1,2}(\Omega):\ 0\leq v\leq 1\text{ a.e. in }\Omega\}.$
We note that the result does not change if in the definition of $P_{\eta}$ we omit the constraint
$$
0\leq v\leq 1\text{ a.e. in }\Omega.
$$
Also the following result, due to Modica, \cite{Mod}, will be useful.

\begin{prop}\label{modcompprop}
For any $C>0$ and any $\eta>0$, let us define
$$A_C=\{v\in L^1(\Omega):\ 0\leq v\leq 1\text{ a.e. and }P_{\eta}(v)\leq C\}.$$
Then $A_C$ is precompact in $L^1(\Omega)$.
\end{prop}

\begin{oss}\label{compactnessoss}
With the same proof, we can show the following. Let us consider
any family $\{v_{\eta}\}_{0<\eta\leq \eta_0}$ such that,
for some positive constant $C$ and
for any $\eta$, $0<\eta\leq\eta_0$,
we have $0\leq v_{\eta}\leq 1$ almost everywhere and
$P_{\eta}(v_{\eta})\leq C$.
Then $\{v_{\eta}\}_{0<\eta\leq \eta_0}$
is precompact in $L^1(\Omega)$.
\end{oss}

Let us fix $q$, $1<q<+\infty$.
Let $V:\mathbb{R}\to[0,+\infty)$ be a continuous function such that $V(t)=0$ if and only if $t=1$
and let $c_V=\int_{0}^1\sqrt{V(t)}\rmd t$.
Let $\psi:\mathbb{R}\to \mathbb{R}$ be a continuous non-decreasing function such that $\psi(0)=0$,
$\psi(1)=1$ and $\psi(t)>0$ if $t>0$.
For any $\eta>0$, let us fix $o_{\eta}=o_{\eta}(q)\geq 0$ such that $\lim_{\eta\to 0^+}o_{\eta}/\eta^{q-1}=0$.
Finally, we define $\psi_{\eta}=(1-o_{\eta})\psi+o_{\eta}$.
Provided $o_{\eta}<1$, we have that
$\psi_{\eta}$ is a continuous, non-decreasing function
such that $\psi_{\eta}(0)=o_{\eta}$ and
$\psi_{\eta}(1)=1$.

For instance, we may choose $V(t)=(t-1)^2/4$ for any $t\in \mathbb{R}$, whence $4c_V=1$. About $\psi$, we may take
$\psi(t)=t^{\gamma}$, $\gamma>0$,
if $t\geq 0$, while $\psi(t)=0$ if $t<0$
Alternatively, we may choose $\psi(t)=0$ if $t<0$, $\psi(t)=-2t^3+3t^2$ for any $t\in [0,1]$, and
$\psi(t)=1$ for any $t>1$. We may finally take $o_{\eta}(q)=\eta^q$.

Then, for any $\eta>0$, we define the following functional
$\mathcal{AT}_{\eta}^q:L^1(\Omega)\times L^1(\Omega)\to [0,+\infty]$ by
\begin{multline}\label{AT}
\mathcal{AT}_{\eta}^q(u,v)=\displaystyle{b\int_{\Omega}\psi_{\eta}(v)|\nabla u|^q+\frac{1}{\eta}
\int_{\Omega}V(v)+\eta \int_{\Omega}|\nabla v|^2}\\
\text{if }u\in W^{1,q}(\Omega)\text{ and }
v\in W^{1,2}(\Omega,[0,1]),
\end{multline}
whereas $\mathcal{AT}_{\eta}^q(u,v)=+\infty$ otherwise.
We shall refer to $\mathcal{AT}_{\eta}^q$ as the Ambrosio-Tortorelli functional.

Let us define the following variant of the Mumford-Shah functional.
The main difference is that we allow the exponent $q$ to be different for $2$, requiring only that $1<q<+\infty$.
For reasons which will appear evident soon, we also add a formal variable $v$ and we pick $c=4c_V$.
We define the
functional $\mathcal{MS}^q:L^1(\Omega)\times L^1(\Omega)\to [0,+\infty]$
by
\begin{multline}\label{MS}
\mathcal{MS}^q(u,v)=\displaystyle{b\int_{\Omega}|\nabla u|^q+4c_V\mathcal{H}^{N-1}(J(u))}\\
\text{if }u\in GSBV(\Omega)\text{ and }
v= 1\text{ a.e. in }\Omega,
\end{multline}
whereas $\mathcal{MS}^q(u,v)=+\infty$ otherwise.

The Ambrosio-Tortorelli functional approximates the Mumford-Shah functional, in the sense of $\Gamma$-convergence.
Such an important approximation result is due to Ambrosio and Tortorelli, \cite{Amb e Tor1,Amb e Tor2}, see also
\cite{Bra}.

\begin{teo}\label{Amb-Torprop}
With respect to the metric of $L^1(\Omega)\times L^1(\Omega)$, we have that, as $\eta\to 0^+$,
$\mathcal{AT}_{\eta}^q$ $\Gamma$-converges to $\mathcal{MS}^q$.
\end{teo}

Let us review some regularity results which will be needed in the sequel. Most of these results are
a consequence of a theorem by Meyers, \cite{Mey}, see also \cite{Gal e Mon}, and of standard regularity estimates, and we shall omit the proofs. Let
$\Omega\subset \mathbb{R}^N$ be a bounded open set with
Lipschitz boundary. Let $A=A(x)$, $x\in \Omega$, be an $N\times N$ matrix whose entries are measurable and such that,
for some $0<\lambda<1$, we have
\begin{equation}
\begin{array}{ll}
A(x)\xi\cdot\xi\geq \lambda |\xi|^2 &\text{for any }\xi\in\mathbb{R}^N\text{ and for a.e. }x\in \Omega,\\
\|A\|_{L^{\infty}(\Omega)}\leq \lambda^{-1}.
\end{array}
\end{equation} 
We remark that for any matrix $A$, by $\|A\|$ we denote the norm of the matrix as a linear operator.

Let $f\in L^s(\partial \Omega)$, with $s>1$ if $N=2$ or
$s\geq 2(N-1)/N$ if $N\geq 3$,
be such that $\int_{\partial \Omega}f=0$ and let
$F\in L^p(\Omega,\mathbb{R}^N)$, with $p\geq 2$.
Let us denote $W^{1,2}_{\ast}(\Omega)=\{u\in W^{1,2}(\Omega):\ \int_{\Omega}u=0\}$.
Then, there exists a unique $u\in W^{1,2}_{\ast}(\Omega)$ such that
\begin{equation}
\int_{\Omega}A\nabla u\cdot\nabla v=
\int_{\Omega}F\cdot \nabla v+
\int_{\partial \Omega}fv\quad\text{for any }v\in W^{1,2}(\Omega).
\end{equation}
This is the weak formulation of
$$\left\{
\begin{array}{ll}
\mathrm{div}(A\nabla u)=\mathrm{div}(F) &\text{in }\Omega\\
A\nabla u\cdot \nu=f  &\text{on }\partial\Omega.
\end{array}
\right.$$

The following regularity result holds true.

\begin{prop}\label{regprop1}
Under the previous assumptions, the following regularity properties hold.

First of all, we have, for a constant $C_0$ depending on
$N$, $\lambda$, $p$, $s$ and $\Omega$ only, 
\begin{equation}\label{poincare}
\|u\|_{W^{1,2}(\Omega)}\leq C_0\left(\|F\|_{L^p(\Omega)}+
\|f\|_{L^s(\partial\Omega)}\right).
\end{equation}

If $p>N$ and $s>N-1$,
there exist a constant $C_1>0$ such that
\begin{equation}\label{infty}
\|u\|_{L^{\infty}(\Omega)}\leq C_1\left(\|F\|_{L^p(\Omega)}+\|f\|_{L^{s}(\partial\Omega)}\right).
\end{equation}
Here $C_1$ depends on $N$, $\lambda$, $p$, $s$ and $\Omega$ only.

There exists a constant $Q>2$, depending on $N$, $\lambda$ and $\Omega$ only \textnormal{(}$Q\to 2$ if $\lambda\to 0^+$\textnormal{)}, such that
if  $p$ satisfies $2<p<Q$ and $s\geq (N-1)p/N$, then
\begin{equation}\label{highintuv}
\|\nabla u\|_{L^{p}(\Omega)}\leq C_2\left(
\|F\|_{L^p(\Omega)}+\|f\|_{L^{s}(\partial\Omega)}\right).
\end{equation}
Here $C_2$ depends on $N$, $\lambda$, $p$, $s$ and $\Omega$ only.

We conclude that if $s>N-1$, 
there exists a constant $q(\lambda)>2$, depending on $N$, $\lambda$, $s$ and
$\Omega$ only, such that for any $p$, $2\leq p\leq q(\lambda)$, we have
\begin{equation}
\|\nabla u\|_{L^{p}(\Omega)}\leq C_3\left(
\|F\|_{L^p(\Omega)}+\|f\|_{L^{s}(\partial\Omega)}\right),
\end{equation}
in particular, if $p=q(\lambda)$
\begin{equation}
\|\nabla u\|_{L^{q(\lambda)}(\Omega)}\leq C_4\left(
\|F\|_{L^{q(\lambda)}(\Omega)}+\|f\|_{L^{s}(\partial\Omega)}\right).
\end{equation}
Here  $C_3$ depends on
$N$, $\lambda$, $p$, $s$ and
$\Omega$ only, whereas $C_4$ depends on
$N$, $\lambda$, $s$ and
$\Omega$ only.
\end{prop}

\begin{oss}
Let us observe that $Q$ and $q(\lambda)$ converges to $2$ as $\lambda\to 0^+$,
whereas all the constants $C_0$--$C_4$ might tend to $+\infty$ as  $\lambda\to 0^+$.
Let us also remark that the same kinds of estimates hold true if we replace $W^{1,2}_{\ast}(\Omega)$ with, for instance,
$${W_E^{1,2}(\Omega)}=\left\{u\in W^{1,2}(\Omega):\ \int_Eu=0\right\}$$
where $E$ is a Borel subset
of $\partial \Omega$ with non-empty interior,
clearly with respect to the induced topology of $\partial\Omega$.
In this case, the constants $C_0$--$C_4$ might depend on $E$ as well.
\end{oss}

We conclude this section with the following lemma, in which we state a Caccioppoli inequality.

\begin{lem}\label{caccioppolilemma}
Let us assume that $A=A(x)$, $x\in B_{2R}$, is a symmetric
$N\times N$ matrix whose entries belong to $L^{\infty}(B_{2R})$. We also assume that,
for some constants $0<\lambda<\Lambda$, we have
$$\lambda|\xi|^2\leq A(x)\xi\cdot\xi\leq \Lambda|\xi|^2\quad\text{for every }\xi\in\mathbb{R}^N\text{ and for a.e. }x\in B_{2R}.$$
Let $w\in L^{\infty}(B_{2R})$ be a weight satisfying $0<\varepsilon\leq w\leq 1$ almost everywhere in $B_{2R}$.
If $u$ solves in a weak sense
$$\mathrm{div}(wA\nabla u)=0\quad\text{in }B_{2R},$$
then
\begin{equation}\label{Caccioppoli}
\int_{B_R}w|\nabla u|^2\leq \frac{C}{R^2}\int_{B_{2R}}wu^2 
\end{equation}
where $C$ depends on $\lambda$ and $\Lambda$ only.
\end{lem}

\proof{.}
In order to prove \eqref{Caccioppoli} it is enough to take a cutoff function $\chi$ such that
$\chi\in C^{\infty}_0(B_{2R})$,
$0\leq \chi\leq 1$ on $B_{2R}$, and $\chi\equiv 1$ on $B_R$. We may also assume that for an absolute constant $C$ we have
$|\nabla \chi|\leq C/R$ on $B_{2R}$.
Then we use the test function $u\chi^2$ in the equation and with simple computations we obtain \eqref{Caccioppoli}.\cvd

\section{The direct problem and setting of the inverse problem}\label{dirsec}

Let $\Omega$, $\Omega_1$ and $\tilde{\Omega}_1$ be
three bounded domains contained in $\mathbb{R}^N$, $N\geq 2$,
with Lipschitz boundaries such that $\Omega_1\subset\tilde{\Omega}_1\subset\Omega$
and the following properties are satisfied. First, $\Omega\backslash\overline{\tilde{\Omega}_1}$ is not empty and 
$\mathrm{dist}(\overline{\Omega_1},\partial\tilde{\Omega}_1\cap\Omega)>0$.
Then, there exist $\gamma$ and $\tilde{\gamma}$, closed subsets of $\partial\Omega$, which
are contained in the interior of $\partial\Omega\cap\partial\Omega_1$ and whose interiors,
with respect to the induced topology of $\partial\Omega$, are not empty.

We assume that $\Omega$, $\Omega_1$, $\tilde{\Omega}_1$, $\gamma$ and $\tilde{\gamma}$
are fixed throughout the paper.

Let $K_0$ be an \emph{admissible defect}, that is $K_0$ is a non-empty compact set contained in $\overline{\Omega}$
such that $\mathrm{dist}(K_0,\overline{\tilde{\Omega}_1})>0$.
We denote with $G_{K_0}$ the connected component of
$\Omega\backslash K_0$ such that $\tilde{\Omega}_1\subset G_{K_0}$. We observe that
$\gamma\cup\tilde{\gamma}\subset \partial G_{K_0}$.
We remark that if $K_0\subset\partial \Omega$ then $G_{K_0}=\Omega$,
that is no defect is present in the conductor.

We say that an admissible defect $K_0$ is a \emph{material loss defect}, or \emph{material loss} for short, if $G_{K_0}$ is equal to the interior of its own closure (that is
no crack-type defect is allowed).

\begin{oss}
If the defect $K_0$ to be  reconstructed is compactly contained in $\Omega$ and $\partial\Omega$ is connected, then we may take $\gamma$ and $\tilde{\gamma}$ equal to $\partial \Omega$. Furthermore, if $\partial \Omega$ is regular enough and
we a priori know that $\mathrm{dist}(K_0,\partial\Omega)>\delta$, for some $\delta>0$ small enough,
then we may choose $\Omega_1=\{x\in \Omega:\ \mathrm{dist}(x,\partial\Omega)<\delta/2 \}$
and $\tilde{\Omega}_1=\{x\in \Omega:\ \mathrm{dist}(x,\partial\Omega)<3\delta/4 \}$.
\end{oss}

Let us fix a number $s$, $s>N-1$, which shall be kept fixed thorughout the paper.
Let us prescribe $f_0\in L^s(\partial\Omega)$ such that $\int_{\partial\Omega}f_0=0$, $f_0\not\equiv 0$
and $\mathrm{supp}(f_0)\subset\tilde{\gamma}$.

Let the electrostatic potential $u_0=u(f_0,K_0)$
be the weak solution to the following Neumann boundary value problem
\begin{equation}\label{dirpbm3}
\left\{\begin{array}{ll}
\Delta u_0=0&\text{in }\Omega\backslash K_0,\\
\nabla u_0\cdot \nu=f_0&\text{on }\tilde{\gamma},\\
\nabla u_0\cdot \nu=0&
\text{on }\partial(\Omega\backslash K_0)\backslash\tilde{\gamma},
\end{array}
\right.
\end{equation}
with the normalization conditions
\begin{equation}\label{normcond1}
\int_{\gamma}u_0=0,
\end{equation}
and
\begin{equation}\label{normcond2}
u_0=0\quad\text{almost everywhere in }\Omega\backslash G_{K_0}.
\end{equation}

Let us recall that our measured additional information is $g_0\in L^2(\gamma)$ where $g_0=u_0|_{\gamma}$. By
\eqref{normcond1}, we have $\int_{\gamma}g_0=0$. 

We observe that there exists a unique solution $u_0$
to
\eqref{dirpbm3}-\eqref{normcond1}-\eqref{normcond2} and that it satisfies the
following regularity properties, 
see \cite{Ron07} for further details.

There exists a constant $C_1>0$, depending on $s$, $\Omega$, $\Omega_1$, $\tilde{\Omega}_1$,
$\gamma$ and $\tilde{\gamma}$ only,
such that
\begin{eqnarray}\label{boundnablau}
&&\|\nabla u_0\|_{L^2(\Omega\backslash K_0)}\leq C_1\|f_0\|_{L^s(\tilde{\gamma})},\\
\label{globalboundu}
&&\|u_0\|_{L^{\infty}(\Omega)}\leq C_1\|f_0\|_{L^s(\tilde{\gamma})}.
\end{eqnarray}



The estimate \eqref{globalboundu} guarantees that $u_0$
belongs to $W^{1,2}(\Omega\backslash K_0)$. Furthermore, under the additional assumption that
$\mathcal{H}^{N-1}(K_0)<+\infty$, or equivalently that $\mathcal{H}^{N-1}(\partial G_{K_0})<+\infty$, 
we have that
$u_0$ belongs to $SBV(\Omega)$, 
its approximate discontinuity set
$J(u_0)$ satisfies $\mathcal{H}^{N-1}(J(u_0)\backslash \partial G_{K_0})=0$ and, finally, $\nabla u_0$, the weak derivative of $u_0$ in $\Omega\backslash K_0$,
coincides almost everywhere in $\Omega$ with the approximate gradient
of $u_0$, see for instance
\cite[Proposition~4.4]{Amb e Fus e Pal}.

For any $r$, $1< r<+\infty$, and any Borel set $E\subset\partial\Omega$ whose interior, in the induced topology, is not empty, we define
$${W^{1,r}_E(\Omega)}=\left\{u\in W^{1,r}(\Omega):\ \int_{E}u=0\right\}.$$
We observe that, by a generalized Poincar\'e inequality, on 
${W^{1,r}_{E}(\Omega)}$ the usual $W^{1,r}(\Omega)$ norm and the norm
$\|u\|_{{W^{1,r}_{E}(\Omega)}}=\|\nabla u\|_{L^r(\Omega)}$ are equivalent. Therefore, we shall
set this second one as the natural norm of ${W^{1,r}_{E}(\Omega)}$.

Let us consider a weight $w$ in $\Omega$ satisfying the following properties.
We assume that $w\in L^{\infty}(\Omega)$
and that $w\geq \varepsilon$ almost everywhere in $\Omega$, for some $\varepsilon>0$.

For any such weight $w$, and any $u_1$, $u_2\in W^{1,2}(\Omega)$, we define the bilinear form
$$\langle u_1,u_2\rangle_{w}=\int_{\Omega}w\nabla u_1\cdot\nabla u_2$$
and we denote the seminorm
$$|u_1|_{w}=\langle u_1,u_1\rangle_{w}^{1/2}=
\left(\int_{\Omega}w|\nabla u_1|^2\right)^{1/2}.$$
We denote, for any $u_1\in W^{1,2}(\Omega)$,
$$\|u_1\|_{w}=\left(\|u_1\|_{L^2(\Omega)}^2+|u_1|_{w}^2\right)^{1/2}.$$
We have that
$\|\cdot\|_{w}$ is an equivalent norm for $W^{1,2}(\Omega)$, and
$\langle\cdot,\cdot\rangle_{w}$ is a scalar product on 
${W^{1,2}_{E}(\Omega)}$ whose corresponding norm, 
$|\cdot|_{w}$,
is an equivalent norm for ${W^{1,2}_{E}(\Omega)}$, for any set $E$ as before.

For any such weight $w$ and any $f\in L^s(\partial\Omega)$ such that
$\int_{\partial\Omega}f=0$ and $\mathrm{supp}(f)\subset\tilde{\gamma}$,
let $u=u(w)$ be the solution to the following Neumann type boundary value problem
\begin{equation}\label{weightpbm}
\left\{\begin{array}{ll}
\mathrm{div}(w\nabla u)=0 &\text{in }\Omega\\
w\nabla u\cdot \nu=f &\text{on }\partial\Omega\\
\int_{\gamma}u=0.
\end{array}
\right.
\end{equation}
The weak formulation of \eqref{weightpbm}
is the following. We look for a function $u\in W^{1,2}_{\gamma}(\Omega)$ such that
$$\int_{\Omega} w\nabla u\cdot\nabla u_1=\int_{\tilde{\gamma}}fu_1\quad\text{for any }u_1\in W^{1,2}(\Omega).$$
Obviously we have existence and uniqueness of such a solution.
Furthermore, the following regularity result holds for $u$.

\begin{prop}\label{crucialprop}
Under the previous notation and assumptions, let $u$ solve \eqref{weightpbm}
for some weight $w$. We assume that
$$\|w\|_{L^{\infty}(\Omega)}\leq A$$
and
$$w(x)=1\quad\text{for a.e. }x\in\tilde{\Omega}_1.$$

Then there exists a constant $C_2$, depending on $s$, $\Omega$, $\Omega_1$, $\tilde{\Omega}_1$, $\gamma$,
$\tilde{\gamma}$ and $A$ only,
such that
\begin{align}
\label{L2uw}
&|u|_w\leq C_2\|f\|_{L^s(\tilde{\gamma})}\\
\label{globalbounduw}
&\|u\|_{L^{\infty}(\Omega)}\leq C_2\|f\|_{L^s(\tilde{\gamma})}.
\end{align}
We notice that the constant $C_2$ does not depend on $w$ or on $\varepsilon$.
\end{prop}

\proof{.} We sketch the proof of this proposition. Inequality \eqref{L2uw} follows from
an application of Poincar\'e inequality in $\Omega_1$.
The $L^{\infty}$ bound \eqref{globalbounduw} is a consequence of the maximum principle and may be proved following the same arguments used to prove \eqref{globalboundu}, see \cite{Ron07} for details.\cvd

\bigskip

Let us fix the notation for our inverse problem. Let $K_0$ be the unknown defect, which for the time being we assume to
be just an admissible defect.

We assume that $f_0$ belongs to $L^s(\partial\Omega)$ and satisfies
$\mathrm{supp}(f_0)\subset\tilde{\gamma}$ and $\int_{\partial\Omega}f_0=0$.
We recall that $s$ is a fixed constant such that $s>N-1$. 

The unknown electrostatic potential is
$u_0=u(K_0,f_0)$, solution to \eqref{dirpbm3}-\eqref{normcond1}-\eqref{normcond2}, and the additional measured data
is $g_0=u_0|_{\gamma}$. We observe that $g_0\in L^2(\gamma)$ and
$\int_{\gamma}g_0=0$.

Let us fix a noise level $\varepsilon$, $0<\varepsilon\leq 1$, then
the noisy Cauchy data are given by
$f_{\varepsilon}$ and $g_{\varepsilon}$. Here $f_{\varepsilon}$ belongs to $L^s(\partial\Omega)$ and satisfies
$\mathrm{supp}(f_{\varepsilon})\subset\tilde{\gamma}$ and $\int_{\partial\Omega}f_{\varepsilon}=0$, whereas
$g_{\varepsilon}$ belongs to $L^2(\gamma)$ and satisfies $\int_{\gamma}g_{\varepsilon}=0$.
We assume that
\begin{equation}\label{noiselevel}
\|f_0-f_{\varepsilon}\|_{L^s(\tilde{\gamma})}\leq \varepsilon \quad\text{and}\quad
\|g_0-g_{\varepsilon}\|_{L^2(\gamma)}\leq \varepsilon.
\end{equation}

For any $0<\varepsilon\leq 1$, let
$\eta=\eta(\varepsilon)>0$ and $a_{\varepsilon}>0$ be such that
$\lim_{\varepsilon\to 0^+}\eta(\varepsilon)=0$ and
$\lim_{\varepsilon\to 0^+}a_{\varepsilon}=0$.
Further assumptions on $\eta(\varepsilon)$ and $a_{\varepsilon}$ will be imposed later.

Let us fix a constant $c_1$, $0<c_1<1$. We recall that $\psi: \mathbb{R}\to \mathbb{R}$ is
a continuous, non-decreasing function such that
$\psi(0)=0$, 
$\psi(1)=1$, and $\psi(t)>0$ if $t>0$. In particular $\psi(c_1)>0$.
Provided $o_{\eta}\leq1/2$, we have that
$\psi_{\eta}$ is a continuous, non-decreasing function
such that $\psi_{\eta}(0)=o_{\eta}$ and
$\psi_{\eta}(1)=1$. Furthermore,
$\psi_{\eta}(c_1)\geq \psi(c_1)/2>0$.

In the sequel we shall always assume that
$$0<o_{\eta}\leq 1/2\quad\text{for any }\eta>0.$$
Without loss of generality, we also assume that $\psi$, $W$, and $V$ are bounded
all over $\mathbb{R}$, for instance by a constant $A$.
For any $\eta>0$,
again without loss of generality,
we assume that $\psi$ is such that
$\psi_{\eta}\geq o_{\eta}/2$
all over $\mathbb{R}$.

To any function $\tilde{v}\in L^1(\Omega)$ we associate the function $v=1-\tilde{v}$. We observe that,
provided $0\leq \tilde{v}\leq 1$
almost everywhere in $\Omega$,
we also have $0\leq v\leq 1$ almost everywhere in $\Omega$.

For any $\eta>0$ and for any $\tilde{v}\in L^1(\Omega)$, let $w_{\eta}=w_{\eta}(\tilde{v})=\psi_{\eta}(v)$,
where $v=1-\tilde{v}$. We observe that $w_{\eta}$ is such that
$\|w_{\eta}\|_{L^{\infty}(\Omega)}\leq A+1/2$ and $w_{\eta}\geq o_{\eta}/2$ almost everywhere in
$\Omega$.
Therefore we define, for any $\tilde{v}\in L^1(\Omega)$ and any $\varepsilon$, $0<\varepsilon\leq 1$,
the function $\tilde{u}_{\varepsilon}\in W^{1,2}_{\gamma}(\Omega)$
where  $\tilde{u}_{\varepsilon}=\tilde{u}_{\varepsilon}(\tilde{v})$ is the solution to the
following boundary value problem
\begin{equation}\label{utilde2}
\left\{\begin{array}{ll}
\mathrm{div}(w_{\eta}(\tilde{v})\nabla \tilde{u}_{\varepsilon})=0 &\text{in }\Omega\\
w_{\eta}(\tilde{v})\nabla \tilde{u}_{\varepsilon}\cdot\nu=f_{\varepsilon} &\text{on }\partial\Omega,
\end{array}\right.
\end{equation}
where as usual $\eta=\eta(\varepsilon)$.

We finally fix positive constants $a_1$, $a_2$, $\tilde{q}$, $\tilde{\beta}$, and $c_2$,
$0<c_1<c_2<1$.
We also define the following space $W(\Omega)=\{\tilde{v}\in W^{1,2}(\Omega):\ \tilde{v}=0\text{ a.e. in }\tilde{\Omega}_1\}$. To any $\tilde{v}\in W(\Omega)$ we associate the function $v=1-\tilde{v}$. We remark that
$v\in W^{1,2}(\Omega)$ and $v=1$ almost everywhere in $\tilde{\Omega}_1$.
All these constants and the notation will be kept fixed throughout the paper.

\section{Determination of material losses}\label{cavssec}

In this section, the main of the paper, we shall consider the problem of determining material losses. We begin by defining
suitable classes of material losses.
 
\begin{defin}\label{defadmmatlos}
Let us fix a positive constant $\delta$. We say that $\mathcal{B}$ is an admissible class of
material losses if the following holds. First, any $K\in \mathcal{B}$ is an admissible defect such
that $\mathrm{dist}(K,\overline{\tilde{\Omega}_1})\geq \delta$ and $G_K$ is a domain with Lipschitz boundary.
Second, we assume that, for some constant $C$,
we have $\mathcal{H}^{N-1}(\partial G_K)\leq C$ for any $K\in \mathcal{B}$. Finally, we assume that
the set $\{\overline{G_K}:\ K\in\mathcal{B}\}$ is compact with respect to the Hausdorff distance.
\end{defin}

In the remaining part of this section, let us fix $\mathcal{B}$, an admissible class of material losses in the
sense of Definition~\ref{defadmmatlos}.
We assume that the unknown defect $K_0$ belongs to $\mathcal{B}$.
We observe that, as in Proposition~\ref{regprop1}, we have  there exist a constant $q>2$ and a constant $C>0$, not depending on
$f_0$, such that
$\nabla u_0\in L^q(\Omega,\mathbb{R}^N)$, in particular
$$\|\nabla u_0\|_{L^q(\Omega)}\leq C\|f_0\|_{L^s(\tilde{\gamma})}.$$
Here the constants $q$ and $C$ depend also on $s$ and on $K_0$.
In the sequel of the section, we shall fix $q>2$ as such a constant, which depends on $K_0$,
among other things. We define $q_1=(q-2)/(2q)$ and we observe that $0<q_1<1/2$. 
We also define the following set.
For any positive constant $a$, we say that
$v\in H(a)$ if $v\in W^{1,2}(\Omega,[0,1])$, $v=1$ almost everywhere in $\tilde{\Omega}_1$
and there exists $K\in \mathcal{B}$ such that $v\geq c_2$ almost everywhere in $\Omega\backslash
\overline{B}_a(\overline{\Omega \backslash G_K})$ and $v\leq c_1$ almost everywhere in
$\Omega\backslash \overline{B}_a(\overline{G_K})$.
We observe that, by the compactness of the class $\mathcal{B}$ with respect to the Hausdorff distance,
such a set $H(a)$ is closed with respect to the weak $W^{1,2}(\Omega)$
convergence.

For any $\varepsilon$, $0<\varepsilon\leq 1$, we define
$\tilde{\mathcal{G}}_{\varepsilon}:W(\Omega)\to \mathbb{R}$
as follows. For any $\tilde{v}\in W(\Omega)$, recalling that $v=1-\tilde{v}$, we set
\begin{equation}
\tilde{\mathcal{G}}_{\varepsilon}(\tilde{v})=
\frac{a_2}{\varepsilon^{\tilde{\beta}}}\int_{\gamma}|\tilde{u}_{\varepsilon}-g_{\varepsilon}|^2+\\
\displaystyle{b\int_{\Omega}w_{\eta}(\tilde{v})|\nabla \tilde{u}_{\varepsilon}|^2+\frac{1}{\eta}
\int_{\Omega}W(v)+\eta \int_{\Omega}|\nabla v|^2}.
\end{equation}
Here $\eta=\eta(\varepsilon)$, $o_{\eta}=o_{\eta}(2)$, $w_{\eta}=w_{\eta(\varepsilon)}(\tilde{v})=
\psi_{\eta(\varepsilon)}(v)$ and  
$\tilde{u}_{\varepsilon}=\tilde{u}_{\varepsilon}(\tilde{v})$ is the solution to \eqref{utilde2}.
Here and in the sequel of this section, we may also set the constant $b=0$, that is we may drop the second term of the functional.

Then, for any $\varepsilon$, $0<\varepsilon\leq 1$, we define
$\mathcal{G}_{\varepsilon}:L^1(\Omega)\to [0,+\infty]$ as follows.
For any $\tilde{v}\in L^1(\Omega)$ we set
\begin{equation}\label{Gepsilon}
\mathcal{G}_{\varepsilon}(\tilde{v})=\tilde{\mathcal{G}}_{\varepsilon}(\tilde{v})\quad\text{if }
\tilde{v}\in W(\Omega)\text{ and }v=(1-\tilde{v})\in H(a_{\varepsilon}),
\end{equation}
whereas $\mathcal{G}_{\varepsilon}(\tilde{v})=+\infty$ otherwise.

\begin{teo}\label{convteo2}
Besides the previous notation and assumptions,
let us further assume that the following constants
satisfy $0<\tilde{\beta}\leq \tilde{q}\leq 2$, and that
$$\limsup_{\varepsilon\to 0^+}\frac{\eta(\varepsilon)^{2q_1}}{\varepsilon^{\tilde{q}}}<+\infty,$$
and, finally, that $a_\varepsilon\geq 2\eta(\varepsilon)$.

Let $u_0=u(f_0,K_0)$. For any $\varepsilon$, $0<\varepsilon\leq 1$,
let 
$$m_{\varepsilon}=\inf\{\mathcal{G}_{\varepsilon}(\tilde{v}):\ \tilde{v}\in L^1(\Omega)\}.$$

Then we have that, for some constant $C$, $m_{\varepsilon}\leq C$ for any $\varepsilon$, $0<\varepsilon\leq 1$. Furthermore, if
$\tilde{v}_{\varepsilon}\in L^1(\Omega)$ is such that
$$\mathcal{G}_{\varepsilon}(\tilde{v}_{\varepsilon})\leq C\quad\text{for any }\varepsilon,\ 0<\varepsilon\leq 1,$$
the following holds.
For any $\varepsilon$, let $v_{\varepsilon}=1-\tilde{v}_{\varepsilon}$ and
$\tilde{u}_{\varepsilon}=\tilde{u}_{\varepsilon}(\tilde{v}_{\varepsilon})$.
Then we have that
$\psi_{\eta(\varepsilon)}(v_{\varepsilon})\tilde{u}_{\varepsilon}\to u_0$ strongly in $L^p(\Omega)$ for any $p$,
$1\leq p <+\infty$, and
$\psi_{\eta(\varepsilon)}(v_{\varepsilon})\nabla \tilde{u}_{\varepsilon}$ converges to $\nabla u_0$ weakly in $L^2(\Omega)$.

Furthermore, for any constant $c$, $c_1<c<c_2$, the sets $\overline{\{v_{\varepsilon}>c\}}$ converge, as
$\varepsilon\to 0^+$,
to $\overline{G_{K_0}}$ in the Hausdorff distance.
\end{teo}

\begin{oss}
We remark that the theorem in particular hold for a family $\tilde{v}_{\varepsilon}\in L^1(\Omega)$ of minimizers or quasi-minimizers, that is satisfying
$$\lim_{\varepsilon\to 0^+}(\mathcal{G}_{\varepsilon}(\tilde{v}_{\varepsilon})-m_{\varepsilon})=0.$$
\end{oss}

\proof{.} By Proposition~\ref{crucialprop}, we infer that there exists a constant $C$
such that for any $\varepsilon$, $0<\varepsilon\leq 1$, and 
for any $\tilde{v}\in W(\Omega)$, we have
$$\int_{\Omega}w_{\eta}(\tilde{v})|\nabla \tilde{u}_{\varepsilon}(\tilde{v})|^2\leq C\quad\text{and}\quad
\|\tilde{u}_{\varepsilon}(\tilde{v})\|_{L^{\infty}(\Omega)}\leq C.$$

By a construction pretty similar to the one used in Proposition~4.5 in \cite{Ron08}, we may construct
$\tilde{v}_{\varepsilon}\in W(\Omega)$ and $u_{\varepsilon}\in W^{1,2}(\Omega)$,
for any $\varepsilon$, $0<\varepsilon\leq 1$, such that the following properties hold.
For any $\varepsilon$, $0<\varepsilon\leq 1$, first
$v_{\varepsilon}=(1-\tilde{v}_{\varepsilon})\in H(a_{\varepsilon})$ and
$$\frac{1}{\eta}
\int_{\Omega}W(v_{\varepsilon})+\eta \int_{\Omega}|\nabla v_{\varepsilon}|^2\leq C.$$
Second, on $\gamma$ we have $u_{\varepsilon}|_{\gamma}=u_0|_{\gamma}=g_0$.
Finally,
$$\int_{\Omega}w_{\eta}(\tilde{v}_{\varepsilon})|\nabla(u_{\varepsilon}-\tilde{u}_{\varepsilon}(\tilde{v}_{\varepsilon}))|^2\leq C\varepsilon^{\tilde{q}}.
$$
By Poincar\'e inequality in $\Omega_1$, we conclude that
\begin{equation}\label{boundary}
\int_{\gamma}|\tilde{u}_{\varepsilon}(\tilde{v}_{\varepsilon})-g_{\varepsilon}|^2\leq 2\left(
\int_{\gamma}|\tilde{u}_{\varepsilon}(\tilde{v}_{\varepsilon})-u_{\varepsilon}|^2+\int_{\gamma}|g_0-g_{\varepsilon}|^2\right)\leq
C(\varepsilon^{\tilde{q}}+\varepsilon^2).
\end{equation}
We immediately conclude that for a constant $C$ we have
$m_{\varepsilon}\leq C$ for any $\varepsilon$, $0<\varepsilon\leq 1$.

For any $n\in\mathbb{N}$, let $\varepsilon_n>0$ be such that
$\lim_n\varepsilon_n=0$ and let $\mathcal{G}_n=\mathcal{G}_{\varepsilon_n}$.
Let $\tilde{v}_n$ be such that $\mathcal{G}_n(\tilde{v}_n)\leq C$, for any $n\in\mathbb{N}$.
Then by Remark~\ref{compactnessoss}, we obtain that, up to a subsequence, $v_n$ converges in $L^1(\Omega)$,
and actually in $L^p(\Omega)$ for any $p$, $1\leq p<+\infty$, and almost everywhere in $\Omega$, to a function
$v$. Such a function $v$ is such that $P(v)$ is finite. Furthermore, by the definition of $H(a)$ and the compactness properties
of $\mathcal{B}$, we may also assume that there exists $K\in\mathcal{B}$ such that $v=1$
almost everywhere in $G_K$ and $v=0$ almost everywhere in $\Omega\backslash G_K$. In other words, $v=\chi_{G_K}$.

Let us call $w_n=w_{\eta(\varepsilon_n)}(\tilde{v}_n)$ and $\tilde{u}_n=\tilde{u}_{\varepsilon_n}(\tilde{v}_n)$.
Let us notice that
$\sqrt{w_n}\nabla \tilde{u}_n$ is uniformly bounded in $L^2(\Omega,\mathbb{R}^N)$, therefore, up to a subsequence,
$\sqrt{w_n}\nabla \tilde{u}_n$ converges to $V\in L^2(\Omega,\mathbb{R}^N)$ weakly in $L^2(\Omega,\mathbb{R}^N)$.

Since $v_n\to \chi_{G_K}$ almost everywhere in $\Omega$, we conclude that also $w_n$ and $\sqrt{w_n}$
converge to $\chi_{G_K}$ almost everywhere in $\Omega$ and in $L^p(\Omega)$ 
for any $p$, $1\leq p<+\infty$. 

For any $B_{2R}(y)\subset (\Omega\backslash\overline{G_K})$, we have that $w_n$ converges to zero almost everywhere
in $B_{2R}(y)$. By the uniform $L^{\infty}$ bound on $\tilde{u}_n$ and by the dominated convergence theorem,
we conclude that $\int_{B_{2R}(y)}w_n\tilde{u}_n^2\to 0$ as $n\to \infty$. By the Caccioppoli inequality described in
Lemma~\ref{caccioppolilemma}, \eqref{Caccioppoli}, we conclude that 
$\sqrt{w_n}\nabla \tilde{u}_n$ converges to $0$ strongly in $L^2(B_R(y),\mathbb{R}^N)$, consequently
$V=0$ almost everywhere in $\Omega\backslash\overline{G_K}$. We conclude that
$w_n\nabla \tilde{u}_n$ weakly converges to $V$ in $L^2(\Omega,\mathbb{R}^N)$ as well.
On the other hand, again up to subsequences and by using the property of $H(a)$, we may follow the arguments of the proof of
Proposition~4.3 in \cite{Ron08} in order to
find a function
$\tilde{u}$ with the following properties. First, $\sqrt{w_n}\tilde{u}_n$ converges to $\tilde{u}$
almost everywhere in $\Omega$ and consequently in $L^p(\Omega)$ for any $p$, $1\leq p<+\infty$.
Second, by the same reasoning above, we conclude that $\tilde{u}=0$ almost everywhere in $\Omega\backslash\overline{G_K}$
and that also $w_n\tilde{u}_n$ converges to $\tilde{u}$
almost everywhere in $\Omega$ and in $L^p(\Omega)$ for any $p$, $1\leq p<+\infty$.
Then, we have that $\tilde{u}\in W^{1,2}(G_K)$, $\tilde{u}$ is harmonic in $G_K$
and $\nabla \tilde{u}=V$ in $G_K$. We  also have that on $\gamma$
and on $\tilde{\gamma}$, $\tilde{u}_n$ converges to $\tilde{u}$ strongly in $L^p(\gamma\cap\tilde{\gamma})$
for any $p$, $1\leq p<+\infty$. As a consequence, $\tilde{u}=g_0$ on $\gamma$.

Let us take any function $\varphi\in W^{1,2}(G_K)$. Since $G_K$ is a domain with Lipschitz boundary,
therefore it is an extension domain, we can find a function $\tilde{\varphi}\in W^{1,2}(\Omega)$
such that $\tilde{\varphi}=\varphi$ on $G_K$.
We conclude that for any $n\in\mathbb{N}$ we have
$$\int_{\Omega}w_n\nabla\tilde{u}_n\cdot\nabla \tilde{\varphi}=\int_{\tilde{\gamma}}f_{\varepsilon_n}\tilde{\varphi}.$$
Since, as $n\to\infty$,
$$\int_{\Omega}w_n\nabla\tilde{u}_n\cdot\nabla \tilde{\varphi}\to \int_{G_K}\nabla\tilde{u}\cdot \nabla \tilde{\varphi}$$
and
$$\int_{\tilde{\gamma}}f_{\varepsilon_n}\tilde{\varphi}\to \int_{\tilde{\gamma}}f_0\tilde{\varphi},$$
we conclude that
$$\int_{G_K}\nabla\tilde{u}\cdot \nabla \varphi=\int_{\tilde{\gamma}}f_0\varphi\quad\text{for any }\varphi\in W^{1,2}(G_K).$$
Then $\tilde{u}$ solves \eqref{dirpbm3}-\eqref{normcond1}-\eqref{normcond2} with $K_0$ replaced by $K$.
Then, since $\tilde{u}=u_0$ on $\gamma$, we conclude by using Theorems 3.3 and 3.6 in \cite{Ron06} that
$\tilde{u}=u_0$ almost everywhere in $\Omega$ and that $G_K=G_{K_0}$. The rest of the proof easily follows.\cvd

\bigskip

We conclude this section with the following existence results, which may be easily proved by the direct method.

\begin{prop}\label{existencemin2}
The following problems admit a solution.
\begin{enumerate}[\textnormal{(}i\textnormal{)}]

\item\label{minpbm3}
$\min\tilde{\mathcal{G}}_{\varepsilon}$
on $W(\Omega)$, with constraint $0\leq \tilde{v}\leq 1$.

\item\label{minpbm4}
$\min\tilde{\mathcal{G}}_{\varepsilon}$
on $W(\Omega)$, with constraints $0\leq \tilde{v}\leq 1$ and
$v\in H(a_{\varepsilon})$ 
\textnormal{(}that is there exists the minimum of $\mathcal{G}_{\varepsilon}$ over
$L^1(\Omega)$\textnormal{)}.
\end{enumerate}
\end{prop}

\section{The crack case}\label{compsec}

In this section we shall deal with the determination of general defects, in particular of cracks.
We begin by recalling results proved in \cite{Ron08}. We include them here for the convenience of the reader and
to compare them with the new results devoted to the determination of material losses, in particular of cavities, which we treated in
Section~\ref{cavssec}.
For what concerns the classes of admissible defects we shall use in this section, 
let us begin with the following definition. We limit ourselves to the two or three-dimensional case, however it is not difficult to see how these  definitions can be generalized to higher dimensions.

If $N=2$, fixed a positive constant $L\geq 1$,
we say that $\Gamma$ is an $L$-\emph{Lipschitz}, or $L$-$C^{0,1}$, \emph{arc} if, up to a rigid transformation,
$\Gamma=\{(x,y)\in\mathbb{R}^2:\ -a/2\leq x\leq a/2,\ y=\varphi_1(x)\}$, where
$L^{-1}\leq a\leq L$  
and $\varphi_1:\mathbb{R}\to \mathbb{R}$ is a Lipschitz map with Lipschitz constant bounded by $L$ and
such that $\varphi_1(0)=0$. For any $\alpha$, $0\leq\alpha\leq 1$, we say that $\Gamma$ is an
$L$-$C^{1,\alpha}$ \emph{arc} if $\varphi_1$ is $C^{1,\alpha}$ and
its $C^{1,\alpha}$ norm is bounded by $L$. The points $(a/2,\varphi_1(a/2))$ and
$(-a/2,\varphi_1(-a/2))$ will be called the \emph{vertices} or \emph{endpoints} of the arc $\Gamma$.

Let us consider now the case $N=3$.
Let $T$ be the closed equilateral
triangle which is contained in the plane $\Pi=\{(x,y,z)\in\mathbb{R}^3:\ z=0\}$
with vertices $V_1=(0,1,0)$, $V_2=(-\sqrt{3}/2,-1/2,0)$ and $V_3=(\sqrt{3}/2,-1/2,0)$ and $T'\subset\mathbb{R}^2$ be its projection on the
plane $\Pi$.
Fixed a positive constant $L\geq 1$,
we
call an $L$-\emph{Lipschitz}, or $L$-$C^{0,1}$, \emph{generalized triangle} a set $\Gamma$ such that, up to a rigid transformation,
$\Gamma=\{(x,y,z)\in\mathbb{R}^3:\ (x,y)\in\varphi(T'),\ z=\varphi_1(x,y)\}$, where 
$\varphi:\mathbb{R}^2\to\mathbb{R}^2$ is a bi-Lipschitz function with constant $L$ such that $\varphi(0)=0$
and $\varphi_1:\mathbb{R}^2\to \mathbb{R}$ is a Lipschitz map with Lipschitz constant bounded by $L$ and
such that $\varphi_1(0)=0$. For any $\alpha$, $0\leq \alpha \leq1$, we say that $\Gamma$ is an
$L$-$C^{1,\alpha}$ \emph{generalized triangle} if $\varphi_1$ is $C^{1,\alpha}$ and
its $C^{1,\alpha}$ norm is bounded by $L$.

In both cases,
the image through $\varphi$ of any vertex or side of $T'$ will be called a generalized vertex or generalized side of $\varphi(T')$,
respectively.
The image on the graph of $\varphi_1$ of one of the generalized vertices of $\varphi(T')$ will be called a \emph{generalized vertex} of $\Gamma$, whereas the
image of one of the generalized sides of $\varphi(T')$ will be called a \emph{generalized side} of $\Gamma$.
We also remark that there exists a constant $L_1>0$, depending on $L$ only, such that we can find $\varphi_2:\mathbb{R}^3\to \mathbb{R}^3$,
a bi-Lipschitz
function with constant $L_1$, such that $\Gamma=\varphi_2(T)$. 

\begin{defin}\label{classesdef}
Let us assume that $\Omega\subset B_R\subset \mathbb{R}^N$, with $R\geq 1$ and $N=2,3$.
For any positive constants $L\geq 1$, $\delta$ and $c$, $c<1$, any $k=0,1$ and $\alpha$, $0\leq\alpha\leq 1$,
such that $k+\alpha\geq 1$, we define
$\mathcal{B}(N,(k,\alpha),L,\delta,c)$ in the following way. We say that
$A\in \mathcal{B}(N,(k,\alpha),L,\delta,c)$ if and only if
$A\subset \overline{B}_{2R}$, there exists a positive integer $n$, depending on $A$, such that
$A=\bigcup_{i=1}^n\Gamma_i$, $\Gamma_i$ an
$L$-$C^{k,\alpha}$ arc (if $N=2$) or generalized triangle (if $N=3$) for any $i=1,\ldots,n$, such that
the following conditions are satisfied:
\begin{enumerate}[i)]
\item for any $i,j\in\{1,\ldots,n\}$ with $i\neq j$, we have that either
$\Gamma_i\cap\Gamma_j$ is not empty or $\mathrm{dist}(\Gamma_i,\Gamma_j)\geq\delta$;
\item for any $i,j\in\{1,\ldots,n\}$ with $i\neq j$, if $\Gamma_i\cap\Gamma_j$ is not empty
then $\Gamma_i\cap\Gamma_j$ is a common endpoint $V$ if $N=2$ and either a common generalized vertex $V$ or a common
generalized side $\gamma$ if $N=3$.
Furthermore, in such a case, for any $x\in \Gamma_i$ we have $\mathrm{dist}(x,\Gamma_j)\geq c|x-V|$
or $\mathrm{dist}(x,\Gamma_j)\geq c\mathrm{dist}(x,\gamma)$, respectively.
\end{enumerate}
\end{defin}

Let us remark that there exists an integer $M$, depending on $N$, $R$, $L$, $\delta$ and $c$ only, such that 
for any $A\in \mathcal{B}(N,(k,\alpha),L,\delta,c)$ we have that
$n\leq M$.

More importantly, we have that any of the classes $\mathcal{B}$
described in Definition~\ref{classesdef}
is non-empty, is composed of non-empty compact sets and it is
compact with respect to the Hausdorff distance. Finally,
if $A$ belongs to any of these classes, then $\mathcal{H}^{N-1}(A)$ is bounded by a constant depending on the class only.

For the time being, let us fix $\mathcal{B}$ as one of the classes of Definition~\ref{classesdef}. We call the constant $k$, $\alpha$, $L$, $\delta$ and $c$ the \emph{a priori data}
related to $\mathcal{B}$.
For any such class $\mathcal{B}$
we call $\mathcal{B}'$ the class of admissible defects $K$ such that
$\mathrm{dist}(K,\overline{\tilde{\Omega}_1})\geq \delta$,
$\mathcal{H}^{N-2}(K\cap \partial\Omega)<+\infty$ and
there exists $A\in \mathcal{B}$ such that $K\subset A$ and
$\mathcal{H}^{N-2}(K\cap \overline{A\backslash K})<+\infty$.

Moreover, we say that $K\in\mathcal{B}'$ satisfies Assumption~A if the following holds.

\begin{assA} We assume that, for any $x_0\in K\cap\Omega$,
there exists $r>0$, depending on $x_0$, such that for any $U$ connected
component of $(\Omega\backslash K)\cap B_r(x_0)$ we can find $r_1>0$, an open set $U_1$, such that
$U\cap B_{r_1}(x_0)\subset U_1\subset U$, and a bijective map $T:U_1\to (-1,1)^N$ such that the following
properties hold. The maps $T$ and $T^{-1}$ are locally Lipschitz and there exists a constant $C$ such that
$\|D T\|$ and $\|DT^{-1}\|$ are bounded by $C$ almost everywhere. By the regularity of $Q=(-1,1)^N$, $T^{-1}$
can be actually extended up to the boundary and we have that $T^{-1}:[-1,1]^{N}\to\mathbb{R}^N$ is a Lipschitz map with Lipschitz
constant bounded by $C$. Furthermore, if we set
$\Gamma=[-1,1]^{N-1}\times\{1\}$, we require that $T^{-1}(\Gamma)=\partial U_1\cap K_0$, $T^{-1}(0,\ldots,0,1)=x_0$ and
$T^{-1}(y)\in \Omega\backslash K$ for any $y\in [-1,1]^N\backslash\Gamma$.

We assume that, for any $x_0\in K\cap\partial\Omega$,
there exists $r>0$, depending on $x_0$, such that for any $U$ connected
component of $(\Omega\backslash K)\cap B_r(x_0)$ we can find $r_1>0$, an open set $U_1$, such that
$U\cap B_{r_1}(x_0)\subset U_1\subset U$, and a bijective map
$T:U_1\to (0,1)\times (-1,1)^{N-1}$ such that the following
properties hold. The maps $T$ and $T^{-1}$ are locally Lipschitz and there exists a constant $C$ such that
$\|D T\|$ and $\|DT^{-1}\|$ are bounded by $C$ almost everywhere. By the regularity of
$Q_1=(0,1)\times (-1,1)^{N-1}$, $T^{-1}$
can be actually extended up to the boundary and we have that $T^{-1}:\overline{Q_1}\to\mathbb{R}^N$ is a Lipschitz map with Lipschitz
constant bounded by $C$. Furthermore, if we set
$\Gamma_1=[0,1]\times[-1,1]^{N-2}\times\{1\}$ and $\Gamma_2=\{0\}\times [-1,1]^{N-1}$, we require that
$T^{-1}(\Gamma_1)=\partial U_1\cap K$,
$T^{-1}(\Gamma_2)=\partial U_1\cap \partial\Omega$, $T^{-1}(0,\ldots,0,1)=x_0$
and
$T^{-1}(y)\in \Omega\backslash K$ for any $y\in \overline{Q_1}\backslash(\Gamma_1\cup\Gamma_2)$.
\end{assA}

In the sequel we shall fix positive constants $L\geq 1$, $\delta$ and $c$, $c<1$, and $\alpha$, $0\leq \alpha\leq 1$.
We also assume that $\Omega\subset B_R$, for some fixed constant $R\geq 1$.
Let $\mathcal{B}=\mathcal{B}(N,(1,\alpha),L,\delta,c)$.
We assume that the unknown defect $K_0$ belongs to $\mathcal{B}'$ and that it satisfies Assumption~A.
We recall that examples of defects satisfying Assumption~A are described in \cite{Ron07,Ron08}.

The next proposition states that the gradient of $u_0$ satisfies a higher integrability property. 

\begin{prop}\label{qprop}
Under the previous assumptions, there exist a constant $q>2$ and a constant $C>0$, which do not depend on $f_0$, such that
$\nabla u_0\in L^q(\Omega,\mathbb{R}^N)$, in particular
$$\|\nabla u_0\|_{L^q(\Omega)}\leq C\|f_0\|_{L^s(\tilde{\gamma})}.$$
\end{prop}

\proof{.} See the proof of Proposition~4.5 in \cite{Ron07}.\cvd

\bigskip

We remark that the constants $q$ and $C$ in Proposition~\ref{qprop} depend also on $s$ and on $K_0$.

For any $a>0$, we call $H_1(a)$ the set of functions $v\in W^{1,2}(\Omega,[0,1])$ such that
$v=1$ almost everywhere in $\tilde{\Omega}_1$
and for some $A\in\mathcal{B}$ we have $v\geq c_1$ almost everywhere in
$\Omega\backslash\overline{B}_a(A)$, where
again $c_1$ is a constant such that $0<c_1<1$.

For any $0<\varepsilon\leq 1$ and any $q\geq 2$, let us define
$\tilde{\mathcal{F}}^q_{\varepsilon}:{W^{1,q}_{\gamma}(\Omega)}\times W(\Omega)\to \mathbb{R}$
as follows. For any $(u,\tilde{v})\in {W^{1,q}_{\gamma}(\Omega)}\times W(\Omega)$, recalling that $v=1-\tilde{v}$, we set
\begin{multline}
\tilde{\mathcal{F}}^q_{\varepsilon}(u,\tilde{v})=
\frac{a_1}{\varepsilon^{\tilde{q}}}|u-\tilde{u}_{\varepsilon}|^2_{w_{\eta}}+
\frac{a_2}{\varepsilon^{\tilde{\beta}}}\int_{\gamma}|u-g_{\varepsilon}|^2+\\
\displaystyle{b\int_{\Omega}\psi_{\eta}(v)|\nabla u|^q+\frac{1}{\eta}
\int_{\Omega}V(v)+\eta \int_{\Omega}|\nabla v|^2}.
\end{multline}
Here $\eta=\eta(\varepsilon)$, $o_{\eta}=o_{\eta}(q)$, $w_{\eta}=w_{\eta(\varepsilon)}(\tilde{v})=
\psi_{\eta(\varepsilon)}(v)$ and  
$\tilde{u}_{\varepsilon}=\tilde{u}_{\varepsilon}(\tilde{v})$ is the solution to \eqref{utilde2}.
We also recall that
$$|u-\tilde{u}_{\varepsilon}|^2_{w_{\eta}}=\int_{\Omega}\psi_{\eta(\varepsilon)}(v)|\nabla (u-\tilde{u}_{\varepsilon})|^2=
\int_{\Omega}\psi_{\eta(\varepsilon)}(v)|\nabla u|^2-2
\int_{\tilde{\gamma}}f_{\varepsilon} u+
\int_{\tilde{\gamma}}f_{\varepsilon} \tilde{u}_{\varepsilon}.$$

Then, for any $0<\varepsilon\leq 1$ and any $q\geq 2$, we define
$\mathcal{F}^q_{\varepsilon}$ as the following functional on
$L^1(\Omega)\times L^1(\Omega)$.
For any $(u,\tilde{v})\in L^1(\Omega)\times L^1(\Omega)$ we set
\begin{multline}\label{Fepsilon}
\mathcal{F}^q_{\varepsilon}(u,\tilde{v})=\tilde{\mathcal{F}}^q_{\varepsilon}(u,\tilde{v})\\ \text{if }
(u,\tilde{v})\in W^{1,q}_{\gamma}(\Omega)\times W(\Omega)\text{ and }v=(1-\tilde{v})\in H_1(a_{\varepsilon}),
\end{multline}
whereas $\mathcal{F}^q_{\varepsilon}(u,\tilde{v})=+\infty$ otherwise.


Now we shall fix the constant $q>2$ as the one defined in Proposition~\ref{qprop}, which depends on $K_0$, among other things. 
Again we set $q_1=(q-2)/(2q)$ and we observe that $0<q_1<1/2$.
The following convergence result is the main result of \cite{Ron08}.

\begin{teo}\label{convteo}
Besides the previous notation and assumptions,
let us further assume that the following constants
satisfy $0<\tilde{q}\leq 2$, $0<\tilde{\beta}\leq 2$, and that
$$\limsup_{\varepsilon\to 0^+}\frac{\eta(\varepsilon)^{2q_1}}{\varepsilon^{\tilde{q}}}<+\infty,$$
and, finally, that $a_\varepsilon\geq 2\eta(\varepsilon)$.

Let $u_0=u(f_0,K_0)$. Then there exists a constant $E_0$, $E_0$ depending on $s$, $\Omega$, $\Omega_1$, $\tilde{\Omega}_1$,
$\gamma$ and $\tilde{\gamma}$ only, such that for any $E$, $E_0\leq E<+\infty$, the following holds.

For any $0<\varepsilon\leq 1$, let
$$m_{\varepsilon}=\inf\{\mathcal{F}^q_{\varepsilon}(u,\tilde{v}):\ (u,\tilde{v})\in L^1(\Omega)\times L^1(\Omega)\text{ and }\|u\|_{L^{\infty}(\Omega)}\leq E\}.$$
Then we have that, for some constant $C$, $m_{\varepsilon}\leq C$ for any $0<\varepsilon\leq 1$.

For any $n\in\mathbb{N}$, let $\varepsilon_n>0$ be such that
$\lim_n\varepsilon_n=0$ and
let $(u_n,\tilde{v}_n)\in L^1(\Omega)\times L^1(\Omega)$ be such that
$\|u_n\|_{L^{\infty}(\Omega)}\leq E$ and
$$\mathcal{F}^q_{\varepsilon_n}(u_n,\tilde{v}_n)\leq C\quad\text{for any }n\in\mathbb{N}.$$
Then, up to a subsequence,
$u_n\to u$ strongly in $L^p(\Omega)$ for any $p$, $1\leq p <+\infty$, and
$\psi_{\eta(\varepsilon_n)}(v_n)\nabla u_n\to \nabla u$ strongly in $L^p(\Omega)$ for any $2\leq p<q$,
where $u=u_0$ almost everywhere in $G_{K_0}$ and
$\nabla u=\nabla u_0$ almost everywhere in $\Omega$.

Furthermore, there exist compact sets $\tilde{A}\subset\overline{\Omega}$ and $A\in \mathcal{B}$, such that $\tilde{A}\subset A$ and  
$\mathcal{H}^{N-1}(J(u)\backslash \tilde{A})=0$, satisfying the following property. For any constant $c$, $0<c\leq c_1$, the sets $\overline{\{v_n<c\}}$ converge, as $n\to\infty$,
to $\tilde{A}$ in the Hausdorff distance.
\end{teo}

An analogous to Proposition~\ref{existencemin2} holds true, again easily proved by the direct method.

\begin{prop}\label{existencemin}
Let $E_0$ be as in Theorem~\textnormal{\ref{convteo}}.
Then  for any $p$, $2\leq p\leq q$, and any $E$, $E_0\leq E\leq +\infty$, 
the following problems admit a solution.
\begin{enumerate}[\textnormal{(}i\textnormal{)}]

\item\label{minpbm1}
$\min\tilde{\mathcal{F}}^p_{\varepsilon}$
on ${W^{1,p}_{\gamma}(\Omega)}\times W(\Omega)$, with constraints $0\leq \tilde{v}\leq 1$
and $\|u\|_{L^{\infty}(\Omega)}\leq E$.

\item\label{minpbm2}
$\min\tilde{\mathcal{F}}^p_{\varepsilon}$
on ${W^{1,p}_{\gamma}(\Omega)}\times W(\Omega)$, with constraints $0\leq \tilde{v}\leq 1$,
$v\in H_1(a_{\varepsilon})$ and $\|u\|_{L^{\infty}(\Omega)}\leq E$
\textnormal{(}that is there exists the minimum of $\mathcal{F}^p_{\varepsilon}$ over
$L^1(\Omega)\times L^1(\Omega)$ with the same $L^{\infty}$ bound on $u$\textnormal{)}.
\end{enumerate}
\end{prop}

Let us now consider the main differences between the cracks and material losses cases.
Our aim is to show the optimality of Theorem~\ref{convteo}, by showing that a reduction to a functional depending on the phase-variable only, with similar convergence properties, may not be feasible. As we have shown in the previous section such a reduction is instead possible in the material loss case.

By Proposition~\ref{regprop1} and Proposition~\ref{crucialprop}, we infer that there exists a constant $C$
such that for any $\varepsilon$, $0<\varepsilon\leq 1$, and 
for any $\tilde{v}\in W(\Omega)$, we have
$$\int_{\Omega}w_{\eta}(\tilde{v})|\nabla \tilde{u}_{\varepsilon}(\tilde{v})|^2\leq C\quad\text{and}\quad
\|\tilde{u}_{\varepsilon}(\tilde{v})\|_{L^{\infty}(\Omega)}\leq C.$$
Furthermore, there exists $q(\varepsilon)>2$, depending on $N$, $\Omega$, $s$ and $\varepsilon$ only, such that
$\tilde{u}_{\varepsilon}(\tilde{v})$
belongs to $W^{1,q(\varepsilon)}(\Omega)$.
We can also find a constant $C_1$, depending on $N$, $\Omega$, $\gamma$, $s$, $\|f_0\|_{L^s(\partial\Omega)}$ and $\varepsilon$ only, such that for any $\tilde{v}\in W(\Omega)$
$$\|\nabla \tilde{u}_{\varepsilon}(\tilde{v})\|_{L^{q(\varepsilon)}(\Omega)}\leq C_1.$$
We remark that the dependence of $q(\varepsilon)$ on $\varepsilon$ is through $o_{\eta(\varepsilon)}$ and that, unfortunately, it might happen that $q(\varepsilon)\to 2^+$
and $C_1\to +\infty$ as $\varepsilon\to 0^+$.

Let us consider the following operator. For any $\varepsilon$, $0<\varepsilon\leq 1$,
we define $\mathcal{H}_{\varepsilon}:W(\Omega)\to W^{1,2}_{\gamma}(\Omega)$ as follows
$$\mathcal{H}_{\varepsilon}(\tilde{v})=\tilde{u}_{\varepsilon}(\tilde{v})\quad\text{for any }\tilde{v}\in W(\Omega).$$
We recall that for any $r$, $1<r<+\infty$, we endow $W^{1,r}_{\gamma}(\Omega)$
with the norm $\|u\|_{W^{1,r}_{\gamma}(\Omega)}=\|\nabla u\|_{L^r(\Omega)}$ for any $u\in W^{1,r}_{\gamma}(\Omega)$. We observe that $\mathcal{H}_{\varepsilon}$
is continuous with respect to the weak-$W^{1,2}(\Omega)$ convergence in $W(\Omega)$ and strong convergence in $W^{1,2}_{\gamma}(\Omega)$.

We obtain that for any $q$, $2\leq q\leq q(\varepsilon)$, we have that
$\mathcal{H}_{\varepsilon}:W(\Omega)\to {W^{1,q}_{\gamma}(\Omega)}$
and that for any $q$, $2\leq q< q(\varepsilon)$,
$\mathcal{H}_{\varepsilon}$
is continuous again with respect to the weak-$W^{1,2}(\Omega)$ convergence in $W(\Omega)$  and strong convergence in $W^{1,q}_{\gamma}(\Omega)$.

Then for any $q\geq 2$, let us define
$\hat{\mathcal{F}}^q_{\varepsilon}:W(\Omega)\to [0,+\infty]$ as follows. For any $\tilde{v}\in W(\Omega)$ we set
\begin{multline}\label{hatdef}
\hat{\mathcal{F}}^q_{\varepsilon}(\tilde{v})=\tilde{\mathcal{F}}^q_{\varepsilon}(\mathcal{H}_{\varepsilon}(\tilde{v}),\tilde{v})=\\
\frac{a_2}{\varepsilon^{\tilde{\beta}}}\int_{\gamma}|\tilde{u}_{\varepsilon}(\tilde{v})-g_{\varepsilon}|^2+
\displaystyle{b\int_{\Omega}\psi_{\eta}(v)|\nabla \tilde{u}_{\varepsilon}(\tilde{v})|^q+\frac{1}{\eta}
\int_{\Omega}V(v)+\eta \int_{\Omega}|\nabla v|^2}.
\end{multline}
Let us notice that for any $q\geq 2$, we have that
there exists $\min \hat{\mathcal{F}}^q_{\varepsilon}$ on $W(\Omega)$
with the constraint $0\leq \tilde{v}\leq 1$,
and with the constraints $0\leq \tilde{v}\leq 1$ and $v\in H_1(a_{\varepsilon})$ as well. 

We investigate whether, for some $q\geq 2$, we may have convergence properties for $\hat{\mathcal{F}}^q_{\varepsilon}$ as we have for $\tilde{\mathcal{F}}^q_{\varepsilon}$.
We observe that $\tilde{\mathcal{G}}_{\varepsilon}$ is equal to $\hat{\mathcal{F}}^2_{\varepsilon}$ but to replace the single-well potential $V$with the double-well potential $W$. It would be desirable to have a convergence result for $\hat{\mathcal{F}}^2_{\varepsilon}$, or at least for $\hat{\mathcal{F}}^q_{\varepsilon}$ with some $q>2$, as we have
for $\tilde{\mathcal{G}}_{\varepsilon}$, Theorem~\ref{convteo2}.
By counterexamples we show that difficulties arise in both cases.
We begin with the case $q=2$ and then we deal with the case $q>2$.

By the construction used in \cite[Proposition~4.5]{Ron08} the next proposition immediately follows.

\begin{prop}\label{recoveryprop} Under the
assumptions of Theorem~\textnormal{\ref{convteo}} and if $0<\tilde{\beta}\leq \tilde{q}\leq 2$,
we can find $\tilde{v}_{\varepsilon}$ for any $\varepsilon$, $0<\varepsilon\leq1$, such that the following holds. For any $\varepsilon$, $0<\varepsilon\leq1$, we have, first, that
$$\hat{\mathcal{F}}^2_{\varepsilon}(\tilde{v}_{\varepsilon})\leq C.$$
Second, $\{v_{\varepsilon}<1/2\}=\{x\in\Omega:\ \mathrm{dist}(x,K_0)<\xi_{\eta}+\eta/2\}$
where $\xi_{\eta}=\sqrt{\eta o_{\eta}}$. 
Finally,
for any $n\in\mathbb{N}$, let $\varepsilon_n>0$ be such that
$\lim_n\varepsilon_n=0$ and let $\tilde{v}_{n}=\tilde{v}_{\varepsilon_n}$ and
$\tilde{u}_n=\tilde{u}_{\varepsilon_n}(\tilde{v}_n)$.
Then, up to a subsequence,
$\tilde{u}_n\to u$ strongly in $L^p(\Omega)$ for any $p$, $1\leq p <+\infty$, and
$\psi_{\eta(\varepsilon_n)}(v_n)\nabla \tilde{u}_n\to \nabla u$ strongly in $L^2(\Omega)$, where $u=u_0$ almost everywhere in $G_{K_0}$ and
$\nabla u=\nabla u_0$ almost everywhere in $\Omega$.
\end{prop}

In terms of $\Gamma$-convergence, we have obtained a kind of $\Gamma$-limsup inequality. What is missing is the corresponding $\Gamma$-liminf inequality, because taking $q=2$ does not guarantee enough compactness. In fact the solutions to the corresponding weighted elliptic problems may converge to a function which is not a solution to a material loss direct problem, as we shall show in Example~\ref{example1} where we use the
instability of the Neumann problem with respect to boundary variations.

In any case, trying to solve the inverse problem by minimizing
$\hat{\mathcal{F}}^2_{\varepsilon}$ on $W(\Omega)$
with the constraints $0\leq \tilde{v}\leq 1$ and $v\in H_1(a_{\varepsilon})$, might be a good strategy. We recall that in this case the assumption $0<\tilde{\beta}\leq \tilde{q}\leq 2$ should be adopted. In fact, minimizing
$\hat{\mathcal{F}}^2_{\varepsilon}$ is numerically simpler than minimizing
$\tilde{\mathcal{F}}^q_{\varepsilon}$ and  still leads to good numerical reconstructions. In fact this method is adopted in \cite{Rin e Ron} and the numerical simulations presented there show its efficacy. 

\begin{exam}\label{example1}
Let us consider the following example. Let $D$ be a smooth bounded domain of $\mathbb{R}^{N-1}$, $N\geq 2$, and let $\lambda^2>0$ be a Neumann eigenvalue for $-\Delta$ on $D$ and let $f$ be a corresponding eigenfunction, that is
$$\left\{\begin{array}{ll}
-\Delta f=\lambda^2 f &\text{in }D\\
\nabla f\cdot \nu=0 &\text{on }\partial D.
\end{array}\right.
$$
We notice that $\int_{D}f=0$ and we may normalize $f$ in such a way that $\int_{D}|f|^2=1$.

For some constant $T>2$, to be fixed later, let $\Omega=D\times (0,T)$ and let
$G_{K_0}=D\times (0,2)$, that is $K_0=\overline{D}\times \{2\}$. Let $\gamma=\tilde{\gamma}=\overline{D}\times\{0\}$.
Then let $u_0$ be a solution to
$$\left\{\begin{array}{ll}
\Delta u_0=0 &\text{in }\Omega\\
\nabla u_0\cdot \nu=f &\text{on }D\times \{0\}\\
\nabla u_0\cdot \nu=0 &\text{on }K_0\\
\nabla u_0\cdot \nu=0 &\text{on }\partial D\times(0,T).
\end{array}\right.
$$
We normalize $u_0$ in such a way that $\int_{\gamma}u_0=0$ and,
by separation of variables, we have that
$$u_0(x,y)=\frac{f(x)}{\lambda}\left[\frac{\cosh(2\lambda)}{\sinh(2\lambda)}\cosh(\lambda y)-\sinh(\lambda y)\right],
\quad x\in D,\ y\in (0,2),$$
whereas $u_0$ may be chosen identically equal to $0$ in $D\times (2,T)$.

By a simple computation, again by separation of variables, we may find $T>2$ and $\mu>0$ and two functions $u^-$ and $u^+$ such that the following conditions hold. First, $u^-=u_0$ in $D\times (0,1)$ and $u^+$ solves
$$\left\{\begin{array}{ll}
\Delta u^+=0 &\text{in }D\times (1,T)\\
\nabla u^+\cdot \nu=0 &\text{on }D\times \{T\}\\
\nabla u^+\cdot \nu=0 &\text{on }\partial D\times(1,T).
\end{array}\right.
$$
Second, the following transmission condition holds true on $D\times \{1\}$
$$u_y^-(x,1)=u_y^+(x,1)=\mu(u^+(x,1)-u^-(x,1)),\quad x\in D.$$
By following \cite{Mur}, we may then construct a Neumann sieve $K_{\delta}\subset K_0$, $\delta>0$, such that if $u_{\delta}$ solves
$$\left\{\begin{array}{ll}
\Delta u_{\delta}=0 &\text{in }\Omega\backslash K_{\delta}\\
\nabla u_{\delta}\cdot \nu=f &\text{on }\gamma\\
\nabla u_{\delta}\cdot \nu=0 &\text{on }\partial(\Omega\backslash K_{\delta})\backslash \gamma\\
\int_{\gamma}u_{\delta}=0,
\end{array}\right.
$$
the following holds. We have that, as $\delta\to 0^+$, $u_{\delta}$ converges to $u^-=u_0$ weakly in $H^1(D\times (0,1))$, and strongly in $L^2(D\times(0,1))$,
and $u_{\delta}$ converges to $u^+$ weakly in $H^1(D\times (1,T))$, and strongly in $L^2(D\times(1,T))$ . Therefore, the Cauchy data of $u_{\delta}$ on $\gamma$ converges, for instance in $L^2(\gamma)$, to the Cauchy data of $u_0$ on
$\gamma$. 

By using Proposition~\ref{recoveryprop} to approximate $K_{\delta}$ and $u_{\delta}$, for any $n\in\mathbb{N}$ we can find $\varepsilon_n>0$, $\eta_n>0$ and $\tilde{v}_n$ such that
$$\hat{\mathcal{F}}^2_{\varepsilon_n}(\tilde{v}_n)\leq C\quad\text{for any }n\in\mathbb{N},$$
and that, as $n\to \infty$, the following holds. 
First,
$\varepsilon_n\to 0^+$ and $\eta_n\to 0^+$. Second, if  $\tilde{u}_n=\tilde{u}_{\varepsilon_n}(\tilde{v}_n)$, $n\in\mathbb{N}$, then we have that
$\tilde{u}_n$ converges to $u^-=u_0$ strongly in $L^2(D\times (0,1))$ and
$\tilde{u}_n$ converges to $u^+$ strongly in $L^2(D\times (1,T))$. Furthermore,
$\nabla \tilde{u}_n\cdot \nu|_{\gamma}=\nabla u_0\cdot \nu|_{\gamma}$ for any $n\in\mathbb{N}$ and
$$\|\tilde{u}_n-u_0\|_{L^2(\gamma)}\to 0\quad\text{as }n\to \infty.$$

Therefore, even if the Cauchy data of $u_0$ on $\gamma$ are well approximated by those of 
$\tilde{u}_{\varepsilon_n}(\tilde{v}_n)$, we have that
$v_n=1-\tilde{v}_n$ is small in a region close to the corresponding Neumann sieve which is far away from the actual location of the looked-for defect $K_0$. This example shows also the difficulty in proving a convergence result without imposing
any further condition on the region where $v$ is small.
\end{exam}

On the other hand, one might try to minimize 
$\hat{\mathcal{F}}^q_{\varepsilon}$ on $W(\Omega)$ for some $q>2$.
If we take $q>2$, then compactness and convergence would follow as a simple consequence of Theorem~\ref{convteo}, but we may not guarantee that we can find a sequence
of phase-field functions $\tilde{v}_n$ such that
$\hat{\mathcal{F}}^q_{\varepsilon_n}(\tilde{v}_n)$ is uniformly bounded.

Again we use the constraints $0\leq \tilde{v}\leq 1$ and
$v\in H_1(a_{\varepsilon})$. If one would be able to find $\tilde{v}_{\varepsilon}$,
$0<\varepsilon\leq 1$, such that 
$\hat{\mathcal{F}}^q_{\varepsilon}(\tilde{v}_{\varepsilon})\leq C$
for any $ 0<\varepsilon\leq 1$ for some constant $C$, then by Proposition~4.3 in \cite{Ron08}, we would obtain the results of Theorem~\ref{convteo}, replacing
$\tilde{\mathcal{F}}^q_{\varepsilon}$ with $\hat{\mathcal{F}}^q_{\varepsilon}$,
even allowing $E$ to be equal to $+\infty$.

We  believe that constructing such functions $\tilde{v}_{\varepsilon}$ for some $q>2$ is a difficult task and that minimizing
$\hat{\mathcal{F}}^q_{\varepsilon}$ for some $q>2$ might
lead to a not correct reconstruction. In Proposition~\ref{example2} below we show the difficulty of obtaining such a uniform bound. 

In order to have higher integrability of the gradient of $\tilde{u}_{\varepsilon}(\tilde{v})$,
we need to guarantee that $w_{\eta}(\tilde{v})=\psi_{\eta}(v)$ is a weight satisfying certain properties, for instance those described by Stredulinsky in \cite{Stre}.
An important class of weights for which these properties are satisfied is the so-called
Muckenhoupt class $A_2$.

We recall that 
$w$, a non-negative measurable function over $\mathbb{R}^N$,
is a \emph{weight} if $0<w<+\infty$ almost everywhere and $w$ is locally integrable.
We say that a weight $w$ belongs to the \emph{Muckenhoupt
class} $A_2$ if there exists a constant $C$ such that for any ball $B\subset \mathbb{R}^N$ we have
\begin{equation}\label{A2}
\left(\frac{1}{|B|}\int_Bw\right)\left(\frac{1}{|B|}\int_Bw^{-1}\right)\leq C.
\end{equation}
The best constant $C$ for which \eqref{A2} holds is usually referred to as the $A_2$-constant of $w$. We observe that the $A_2$-constant of $w$ is always greater than or equal to $1$. For more details about the Muckenhoupt
weights and weighted elliptic equations, we refer for instance to \cite{Hei e Kil e Mar}.


Therefore, a reasonable assumption is to take $w=w_{\eta}(\tilde{v})$ belonging to
the Muckenhoupt class $A_2$ and such that its $A_2$-constant is bounded by $C$, for some fixed $C$. Without loss of generality we can assume that $0\leq w\leq 1$ almost everywhere and that
$w=1$ outside a given ball $B_{2R}$.
Consequently, we infer that there exists a constant $C_1$, depending on $C$ and
$R$ only, such that
$$\int_{B_R}w^{-1}\leq C_1.$$

\begin{prop}\label{example2}
Let us fix $q>2$.
Let 
$\varepsilon_n$, $n\in\mathbb{N}$, be a sequence of positive numbers such that
$\lim_n\varepsilon_n=0$. For any $n\in\mathbb{N}$, let
$\hat{\mathcal{F}}^q_n=\hat{\mathcal{F}}^q_{\varepsilon_n}$
and let $\tilde{v}_n\in W(\Omega)$ be such that the following holds.
For any $n\in\mathbb{N}$, we assume that $0\leq \tilde{v}_n\leq 1$ and
we set $\eta_n=\eta(\varepsilon_n)$, $v_n=(1-\tilde{v}_n)\in H_1(a_{\varepsilon_n})$,
$w_n=w_{\eta}(\tilde{v}_n)$ and $\tilde{u}_n=\tilde{u}_{\varepsilon_n}(\tilde{v}_n)$.
For any $n\in\mathbb{N}$, we assume
\begin{equation}\label{intinverse}
\int_{\Omega}w_n^{-1}\leq C_1
\end{equation}
and
$$\int_{\Omega}w_n|\nabla\tilde{u}_n|^q+\frac{1}{\eta_n}\int_{\Omega}V(v_n)+\eta_n\int_{\Omega}|\nabla v_n|^2\leq C.$$

Let us consider $\tilde{u}$ as the solution to
\begin{equation}\label{nocrack}
\left\{\begin{array}{ll}
\Delta \tilde{u}=0 &\text{in }\Omega\\
\nabla \tilde{u}\cdot\nu=f_0 &\text{on }\partial\Omega.
\end{array}\right.
\end{equation}
We assume that $K_0$ satisfies the assumption of Theorem~\textnormal{\ref{convteo}}
and that $\tilde{u}\neq u_0$, in particular that
$\tilde{u}|_{\gamma}\neq u_0|_{\gamma}=g_0$. We also assume,
for the time being, that $\Omega$ and $f_0$ are regular enough to guarantee that
$\tilde{u}\in L^{\infty}(\Omega)$ and 
$\nabla\tilde{u}\in L^{\infty}(\Omega,\mathbb{R}^N)$. We may also assume that actually $\tilde{v}_n$ provides a good approximation of $K_0$, namely that $\{v_n<1\}\subset B_{a_n}(K_0)$
where $a_n=a_{\varepsilon_n}$.

Then we have that, as $n\to\infty$,
$w_n\nabla \tilde{u}_n$ converges to $\nabla \tilde{u}$ weakly in $L^2(\Omega)$.
Consequently, as $n\to \infty$, we also have that
$\int_{\gamma}|\tilde{u}_n-g_{\varepsilon_n}|^2\to \int_{\gamma}|\tilde{u}-g_0|^2\neq 0$ and
$\hat{\mathcal{F}}^q_{\varepsilon_n}(\tilde{v}_n)\to +\infty$.
\end{prop}

\proof{.}
Let us compute
$$\int_{\Omega}w_n|\nabla \tilde{u}_n-\nabla\tilde{u}|^2=
\int _{\partial\Omega}(f_{\varepsilon_n}-f_0)(\tilde{u}_n-\tilde{u})+
\int_{\Omega}(1-w_n)\nabla \tilde{u}\cdot\nabla (\tilde{u}_n-\tilde{u}).$$
By the uniform bound on $\tilde{u}_n$, we easily obtain that
$$\int _{\partial\Omega}(f_{\varepsilon_n}-f_0)(\tilde{u}_n-\tilde{u})\to 0\quad\text{as }n\to \infty.$$
Let us evaluate the other term. We have
$$\int_{\Omega}(1-w_n)\nabla \tilde{u}\cdot\nabla (\tilde{u}_n-\tilde{u})=
\int_{\Omega}\frac{(1-w_n)}{\sqrt[q]{w_n}}\sqrt[q]{w_n}\nabla \tilde{u}\cdot\nabla (\tilde{u}_n-\tilde{u}).$$
We apply H\"older inequality with coefficients $q$, $p$ and $r$ such that
$q^{-1}+p^{-1}+r^{-1}=1$ and we obtain
\begin{multline*}
\int_{\Omega}(1-w_n)\nabla \tilde{u}\cdot\nabla (\tilde{u}_n-\tilde{u})\leq\\
\left(\int_{\Omega}w_n^{-r/q}\right)^{1/r}
\left(\int_{\Omega}w_n|\nabla (\tilde{u}_n-\tilde{u})|^q\right)^{1/q}
\left(\int_{\Omega}|1-w_n|^p|\nabla \tilde{u}|^{p}\right)^{1/p}.
\end{multline*}
We use our assumptions to infer
$$
\int_{\Omega}(1-w_n)\nabla \tilde{u}\cdot\nabla (\tilde{u}_n-\tilde{u})\leq C
\left(\int_{\Omega}w_n^{-r/q}\right)^{1/r}
\left(\int_{\Omega}|1-w_n|^p\right)^{1/p},
$$
with $C$ independent of $n$. We may choose $r$ such that $0<r/q<1$,
therefore, since $0\leq w_n\leq 1$, we have
$w_n^{-r/q}\leq w_n^{-1}$. Hence, by \eqref{intinverse}, we conclude that
$\int_{\Omega}w_n|\nabla \tilde{u}_n-\nabla\tilde{u}|^2$ goes to zero as $n\to\infty$.
We obtain that 
$\int_{\gamma}|\tilde{u}_n-\tilde{u}|^2$ goes to zero as well. We then apply Theorem~4.4 in \cite{Ron08} and the
proof is concluded.\cvd

\bigskip

Therefore, even if $v_n$ is a good phase-field approximation of $K_0$,
$\tilde{u}_n$ is not a good approximation of $u_0$. In order to have that
$\tilde{u}_n$ approximates $u_0$, we need to require that $v_n$ is very small close to $K_0$, in such a way that violates \eqref{intinverse}. In turn, this might suggest the fact that higher integrability and the correct approximation might in some sense
oppose each other.

Let us conclude by observing that \eqref{intinverse} is a kind of minimal condition to
have $\int_{\Omega}w_n|\nabla\tilde{u}_n|^q$ uniformly bounded. We wish to point out that potential theory for weights
whose inverse is not integrable has been developed, see for instance \cite{DiF-Zam} for the case of weights
$w=\omega^{1-p/N}$, $1<p<N$, where $\omega$ is a so-called strong $A_{\infty}$-weight. Strong $A_{\infty}$-weights have been introduced in \cite{Dav-Sem}. Following \cite{Sem}, an important example of strong $A_{\infty}$-weights
is given by
$$\omega(x)=\min\{1,\mathrm{dist}(x,A)^s\},\quad x\in\mathbb{R}^N,$$
where $s>0$ and $A$ is a suitable compact set. In \cite[Proposition~4.4]{Sem} it is shown that $\omega$ is a 
strong $A_{\infty}$-weight for any $s>0$ provided $A$ is uniformly disconnected. On the other hand, no strong $A_{\infty}$-weight may vanish on a rectifiable curve, therefore this class of weights seems to be not apt to approximate hypersurfaces as we require in our application.

\section{Differentiability of the functionals}\label{diffsec}

In this last section, we investigate the differentiability properties of
$\tilde{\mathcal{F}}^q_{\varepsilon}$, $\tilde{\mathcal{G}}_{\varepsilon}$ and $\hat{\mathcal{F}}^q_{\varepsilon}$,
for a fixed $\varepsilon$, $0<\varepsilon\leq 1$, and any $q\geq 2$.
For this purpose, we further
assume that the functions $\psi$, $V$ and $W$ are actually of class $C^1$ and such that their
derivatives are bounded and uniformly continuous
all over $\mathbb{R}$.

We define the following spaces. For any $p$, $2\leq p\leq +\infty$,
let us call $L_p(\Omega)=\{\tilde{v}\in L^{p}(\Omega):\tilde{v}=0\text{ a.e. in }\tilde{\Omega}_1\}$
and
 $W_p(\Omega)= W^{1,2}(\Omega)\cap L_{p}(\Omega)$,
 with norm $\|\tilde{v}\|_{L_p(\Omega)}=\|\tilde{v}\|_{L^{p}(\Omega)}$ and
 $\|\tilde{v}\|_{W_p(\Omega)}=\|\tilde{v}\|_{L^{p}(\Omega)}+\|\nabla \tilde{v}\|_{L^{2}(\Omega)}$.
To any $\tilde{v}\in L^{2}(\Omega)$ we as usual associate the function $v=1-\tilde{v}$.
If $\tilde{v}$ belongs either to $L_p(\Omega)$ or to $W_p(\Omega)$, then $v\in L^p(\Omega)$,
$v=1$ almost everywhere in $\tilde{\Omega}_1$, and, provided $0\leq \tilde{v}\leq 1$
almost everywhere in $\Omega$,
we also have $0\leq v\leq 1$ almost everywhere in $\Omega$.
We observe that $W_2(\Omega)=W(\Omega)$ as previously defined. We also recall that
$W^{1,q}_{\gamma}(\Omega)$ is equipped with the norm $\|u\|_{W^{1,q}_{\gamma}(\Omega)}=\|\nabla u\|_{L^q(\Omega)}$
for any $u\in W^{1,q}_{\gamma}(\Omega)$.

We recall that for any $\varepsilon$, $0<\varepsilon\leq 1$,
we define $\mathcal{H}_{\varepsilon}:L_2(\Omega)\to W^{1,2}_{\gamma}(\Omega)$ as follows
$$\mathcal{H}_{\varepsilon}(\tilde{v})=\tilde{u}_{\varepsilon}(\tilde{v})\quad\text{for any }\tilde{v}\in L_2(\Omega).$$

It can be  shown that for any $\tilde{v}_0\in L_2(\Omega)$ such an operator
$\mathcal{H}_{\varepsilon}$ is differentiable in $\tilde{v}_0$ with respect to
the $L^{\infty}(\Omega)$ norm.
Let $D\mathcal{H}_{\varepsilon}(\tilde{v}_0):L_{\infty}(\Omega)\to
W^{1,2}_{\gamma}(\Omega)$ be the differential in $\tilde{v}_0$. Then
for any $\tilde{v}$ in $L_{\infty}(\Omega)$ we have
$$D\mathcal{H}_{\varepsilon}(\tilde{v}_0)[\tilde{v}]=U_{\varepsilon}(\tilde{v}_0,\tilde{v})$$
where $U_{\varepsilon}=U_{\varepsilon}(\tilde{v}_0,\tilde{v})\in {W^{1,2}_{\gamma}(\Omega)}$ solves
the following problem
\begin{equation}\label{der}
\left\{\begin{array}{ll}
\mathrm{div}(\psi_{\eta}(v_0)\nabla U_{\varepsilon})=
\mathrm{div}(\psi'_{\eta}(v_0)\tilde{v}\nabla (\mathcal{H}_{\varepsilon}(\tilde{v}_0)))
 &\text{in }\Omega,\\
\psi_{\eta}(v_0)\nabla U_{\varepsilon}\cdot\nu=0 &\text{on }\partial\Omega.
\end{array}\right.
\end{equation}
Here, obviously, $v_0=1-\tilde{v}_0$.

We recall that for any vector valued function $G\in L^2(\Omega,\mathbb{R}^N)$,
$\mathrm{div}(G)$ defines a functional on $W^{1,2}(\Omega)$ in the following way
$$\mathrm{div}(G)[\phi]=-\int_{\Omega}G\cdot \nabla \phi\quad\text{for any }\phi\in W^{1,2}(\Omega).$$
Therefore,
the weak formulation of \eqref{der} is looking for a function
$U_{\varepsilon}\in{W^{1,2}_{\gamma}(\Omega)}$ such that
$$\int_{\Omega}\psi_{\eta}(v_0)\nabla U_{\varepsilon}\cdot\nabla \varphi
=
\int_{\Omega}\psi'_{\eta}(v_0)\tilde{v}\nabla(\mathcal{H}_{\varepsilon}(\tilde{v}_0))
\cdot \nabla \varphi\quad\text{for any }
\varphi\in {W^{1,2}(\Omega)}.$$

Here, and analogously in the sequel, the differentiability has to be understood in the following sense. For any
$\tilde{v}$ in $L_{\infty}(\Omega)$
$$\mathcal{H}_{\varepsilon}(\tilde{v}_0+\tilde{v})=\mathcal{H}_{\varepsilon}(\tilde{v}_0)+
D\mathcal{H}_{\varepsilon}(\tilde{v}_0)[\tilde{v}]+R(\tilde{v})$$
where
$$\lim_{\|\tilde{v}\|_{L^{\infty}(\Omega)}\to 0}\frac{\|R(\tilde{v})\|_{W^{1,2}_{\gamma}(\Omega)}}{\|\tilde{v}\|_{L^{\infty}(\Omega)}}=0.$$

For any $q\geq 2$,
let us consider the functional
$\tilde{\mathcal{F}}^q_{\varepsilon}:{W^{1,q}_{\gamma}(\Omega)}\times W(\Omega)\to \mathbb{R}$. For  any
$(u_0,\tilde{v}_0)\in {W^{1,q}_{\gamma}(\Omega)}\times W(\Omega)$,
$\tilde{\mathcal{F}}^q_{\varepsilon}$ is differentiable in $(u_0,\tilde{v}_0)$,
with respect to the ${W^{1,q}_{\gamma}(\Omega)}\times W_{\infty}(\Omega)$ norm.
Let $D\tilde{\mathcal{F}}^q_{\varepsilon}(u_0,\tilde{v}_0): {W^{1,q}_{\gamma}(\Omega)}\times W_{\infty}(\Omega)\to \mathbb{R}$
be the differential in $(u_0,\tilde{v}_0)$.
Then, for any $(u,\tilde{v})\in {W^{1,q}_{\gamma}(\Omega)}\times W_{\infty}(\Omega)$,
we have
\begin{multline}\label{differentialf}
D\tilde{\mathcal{F}}^q_{\varepsilon}(u_0,\tilde{v}_0)[(u,\tilde{v})]=\\
\frac{a_1}{\varepsilon^{\tilde{q}}}\int_{\Omega}\left(2\psi_{\eta}(v_0)\nabla u_0\cdot\nabla u
-\psi'_{\eta}(v_0)
|\nabla u_0|^2\tilde{v}
\right)+\\
\frac{a_1}{\varepsilon^{\tilde{q}}}\int_{\tilde{\gamma}}\left(f_{\varepsilon}U_{\varepsilon}(\tilde{v}_0,\tilde{v})
-2f_{\varepsilon}u
\right)+
\frac{2a_2}{\varepsilon^{\tilde{\beta}}}\int_{\gamma}(u_0-g_{\varepsilon})u+\\
b \int_{\Omega}\left(q\psi_{\eta}(v_0)|\nabla u_0|^{q-2}\nabla u_0\cdot\nabla u-
\psi'_{\eta}(v_0)|\nabla u_0|^{q}\tilde{v}\right)+\\
\frac{1}{\eta}\int_{\Omega}(-V'(v_0)\tilde{v})+2\eta
\int_{\Omega}\nabla\tilde{v}_0\cdot\nabla \tilde{v}.
\end{multline}

With the same computation, we infer that
the functionals
$\hat{\mathcal{F}}^2_{\varepsilon}:W(\Omega)\to \mathbb{R}$
and $\tilde{\mathcal{G}}_{\varepsilon}:W(\Omega)\to \mathbb{R}$
are differentiable in $\tilde{v}_0$ for  any
$\tilde{v}_0\in W(\Omega)$, with respect to the $W_{\infty}(\Omega)$ norm.
Let $D\hat{\mathcal{F}}^2_{\varepsilon}(\tilde{v}_0): 
W_{\infty}(\Omega)\to \mathbb{R}$ and
$D\tilde{\mathcal{G}}_{\varepsilon}(\tilde{v}_0): 
W_{\infty}(\Omega)\to \mathbb{R}$
be the differentials in $(u_0,\tilde{v}_0)$.
Then, for any $\tilde{v}\in W_{\infty}(\Omega)$,
we have
\begin{multline}\label{differentialhatf2}
D\hat{\mathcal{F}}^2_{\varepsilon}(\tilde{v}_0)[\tilde{v}]=
\frac{2a_2}{\varepsilon^{\tilde{\beta}}}\int_{\gamma}(\mathcal{H}_{\varepsilon}(\tilde{v}_0)-g_{\varepsilon})U_{\varepsilon}(\tilde{v}_0,\tilde{v})+\\
b \int_{\Omega}\left(2\psi_{\eta}(v_0)\nabla \mathcal{H}_{\varepsilon}(\tilde{v}_0)\cdot\nabla U_{\varepsilon}(\tilde{v}_0,\tilde{v})-
\psi'_{\eta}(v_0)|\nabla \mathcal{H}_{\varepsilon}(\tilde{v}_0)|^{2}\tilde{v}\right)+\\
\frac{1}{\eta}\int_{\Omega}(-V'(v_0)\tilde{v})+2\eta
\int_{\Omega}\nabla\tilde{v}_0\cdot\nabla \tilde{v}
\end{multline}
and
\begin{multline}\label{differentialg}
D\tilde{\mathcal{G}}_{\varepsilon}(\tilde{v}_0)[\tilde{v}]=
\frac{2a_2}{\varepsilon^{\tilde{\beta}}}\int_{\gamma}(\mathcal{H}_{\varepsilon}(\tilde{v}_0)-g_{\varepsilon})U_{\varepsilon}(\tilde{v}_0,\tilde{v})+\\
b \int_{\Omega}\left(2\psi_{\eta}(v_0)\nabla \mathcal{H}_{\varepsilon}(\tilde{v}_0)\cdot\nabla U_{\varepsilon}(\tilde{v}_0,\tilde{v})-
\psi'_{\eta}(v_0)|\nabla \mathcal{H}_{\varepsilon}(\tilde{v}_0)|^{2}\tilde{v}\right)+\\
\frac{1}{\eta}\int_{\Omega}(-W'(v_0)\tilde{v})+2\eta
\int_{\Omega}\nabla\tilde{v}_0\cdot\nabla \tilde{v}.
\end{multline}

It might be useful to have differentiability properties with respect to the $W_{p}(\Omega)$ norm, with $p$ finite.
In fact in this case $W_p(\Omega)$ is a strictly convex real reflexive Banach space and this is useful when we need to apply a gradient method in a numerical implementation, see \cite{Rin e Ron} for details on the use of this information.
In order to obtain such differentiability,
let us now assume that $\psi'$, $V'$ and $W'$
are H\"older continuous for some exponent $\alpha$, $0<\alpha\leq 1$, all over $\mathbb{R}$.

We recall again that, by Proposition~\ref{regprop1} and Proposition~\ref{crucialprop},
there exists a constant $C$
such that for any $\varepsilon$, $0<\varepsilon\leq 1$, and 
for any $\tilde{v}\in L_2(\Omega)$, we have
$$\int_{\Omega}w_{\eta}(\tilde{v})|\nabla \tilde{u}_{\varepsilon}(\tilde{v})|^2\leq C\quad\text{and}\quad
\|\tilde{u}_{\varepsilon}(\tilde{v})\|_{L^{\infty}(\Omega)}\leq C.$$
Furthermore, there exists $q(\varepsilon)>2$, depending on $N$, $\Omega$, $s$ and $\varepsilon$ only, such that
$\tilde{u}_{\varepsilon}(\tilde{v})$
belongs to $W^{1,q(\varepsilon)}(\Omega)$.
We can also find a constant $C_1$, depending on $N$, $\Omega$, $\gamma$, $s$, $\|f_0\|_{L^s(\partial\Omega)}$ and $\varepsilon$ only, such that for any $\tilde{v}\in L_2(\Omega)$
$$\|\nabla \tilde{u}_{\varepsilon}(\tilde{v})\|_{L^{q(\varepsilon)}(\Omega)}\leq C_1.$$
We remark that the dependence of $q(\varepsilon)$ on $\varepsilon$ is through $o_{\eta(\varepsilon)}$ and that, unfortunately, it might happen that $q(\varepsilon)\to 2^+$
and $C_1\to +\infty$ as $\varepsilon\to 0^+$.

However, we may conclude that
$\mathcal{H}_{\varepsilon}:L_2(\Omega)\to W^{1,q(\varepsilon)}_{\gamma}(\Omega)$
and its image is bounded in $W^{1,q(\varepsilon)}_{\gamma}(\Omega)$. Furthermore,
again by Proposition~\ref{regprop1}, we infer that for any $\tilde{v}_0\in L_2(\Omega)$ 
we may define as before
$D\mathcal{H}_{\varepsilon}(\tilde{v}_0)$ and prove that 
$D\mathcal{H}_{\varepsilon}(\tilde{v}_0):L_{q(\varepsilon)(q(\varepsilon)+2)/(q(\varepsilon)-2)}(\Omega)\to W^{1,(q(\varepsilon)+2)/2}_{\gamma}(\Omega)$
is a bounded linear operator.

Let $p(\varepsilon)=q(\varepsilon)\frac{q(\varepsilon)+2}{q(\varepsilon)-2}$. Then, straightforward but lengthy computations allow us to show that
for any $\tilde{v}_0\in L_2(\Omega)$, 
$\mathcal{H}_{\varepsilon}$ is differentiable in $\tilde{v}_0$ with respect to the
$L^{p(\varepsilon)}(\Omega)$ and the $W^{1,2}_{\gamma}(\Omega)$ norms.
The differential is still given by \eqref{der}. We immediately infer that
for any $p$, $p\geq p(\varepsilon)$,
$\tilde{\mathcal{G}}_{\varepsilon}$ is differentiable in $\tilde{v}_0$,
for any $\tilde{v}_0\in W(\Omega)$, with respect to the $W_{p}(\Omega)$ norm, with the
differential given by \eqref{differentialg}.

By an interpolation inequality, we may find $q_1(\varepsilon)$,
$2<q_1(\varepsilon)<(q(\varepsilon)+2)/2$, depending on $q(\varepsilon)$ and
$\alpha$ only, such that
for any $q$, $2\leq q\leq q_1(\varepsilon)$, and any
$p$, $p\geq p(\varepsilon)$, we have that,
for any $\tilde{v}_0\in L_2(\Omega)$, 
$\mathcal{H}_{\varepsilon}$ is differentiable in $\tilde{v}_0$ with respect to the
$L^{p}(\Omega)$ and the $W^{1,q}_{\gamma}(\Omega)$ norms.
Obviously the differential is still given by \eqref{der}.

We conclude that for such $q$ and $p$, and any $\tilde{v}\in W(\Omega)$, we have that
$\hat{\mathcal{F}}^q_{\varepsilon}$ is
is differentiable in $\tilde{v}_0$, for any 
$\tilde{v}_0\in W(\Omega)$, with respect to the $W_p(\Omega)$ norm. Its differential is given by the following formula. For any $\tilde{v}\in W_p(\Omega)$
we have
\begin{multline}\label{differentialhatfq}
D\hat{\mathcal{F}}^q_{\varepsilon}(\tilde{v}_0)[\tilde{v}]=
\frac{2a_2}{\varepsilon^{\tilde{\beta}}}\int_{\gamma}(\mathcal{H}_{\varepsilon}(\tilde{v}_0)-g_{\varepsilon})U_{\varepsilon}(\tilde{v}_0,\tilde{v})+\\
b \int_{\Omega}\left(q\psi_{\eta}(v_0)|\nabla \mathcal{H}_{\varepsilon}(\tilde{v}_0)|^{q-2}\nabla \mathcal{H}_{\varepsilon}(\tilde{v}_0)\cdot\nabla U_{\varepsilon}(\tilde{v}_0,\tilde{v})-
\psi'_{\eta}(v_0)|\nabla \mathcal{H}_{\varepsilon}(\tilde{v}_0)|^{q}\tilde{v}\right)+\\
\frac{1}{\eta}\int_{\Omega}(-V'(v_0)\tilde{v})+2\eta
\int_{\Omega}\nabla\tilde{v}_0\cdot\nabla \tilde{v}.
\end{multline}

An important final remark is the following. If $N=2$, then we may actually choose $p(\varepsilon)=2$,
and we observe that $W_2(\Omega)$
is a Hilbert space, with the scalar product $\int_{\Omega}\nabla \tilde{v}_1\cdot\nabla \tilde{v}_2$ for any
$\tilde{v}_1$, $\tilde{v}_2\in W_2(\Omega)$.
If $N>2$, then it might happen that $p(\varepsilon)>2$ and therefore $W_{p(\varepsilon)}(\Omega)$ 
has not a Hilbert space structure anymore. 
However, since $p(\varepsilon)$ is finite, $W_{p(\varepsilon)}(\Omega)$ is still
a strictly convex real reflexive Banach space.

\end{document}